\documentclass[a4paper]{article}

\usepackage{amsthm}
\usepackage[usenames,dvipsnames]{xcolor}

\newtheorem{thm}{Theorem}[section]

\newtheorem{conjecture}[thm]{Conjecture}
\newtheorem{corollary}[thm]{Corollary}
\newtheorem{definition}[thm]{Definition}
\newtheorem{example}[thm]{Example}
\newtheorem{lemma}[thm]{Lemma}
\newtheorem{open}[thm]{Open}
\newtheorem{proposition}[thm]{Proposition}
\newtheorem{theorem}[thm]{Theorem}

\usepackage{algorithm} 
\usepackage{algpseudocode} 
\usepackage{bm} 
\usepackage{diagbox} 
\usepackage{framed} 
\usepackage{mathdots} 
\usepackage{mathtools} 
\usepackage{tikz} 
\usepackage{hyperref}
\usepackage{amsfonts}

\DeclareMathOperator{\fibon}{Fib}
\newcommand{\I}{\mathcal I}
\DeclareMathOperator{\identityperm}{\I}
\newcommand{\Z}{\mathbb Z}

\newcommand{\calS}{\mathcal S}

\begin{document}

\title{Combinatorial and Asymptotic Results on the Neighborhood Grid}

\author{
    Alex McDonough (University Of California, Davis),\\ 
    Ulrich Reitebuch (FU Berlin),\\
    and Martin Skrodzki (TU Delft)
}

\maketitle

\begin{abstract}
	In various application fields, such as fluid-, cell-, or crowd-simulations, spatial data structures are very important.
	They answer nearest neighbor queries which are instrumental in performing necessary computations for, e.g., taking the next time step in the simulation.
	Correspondingly, various such data structures have been developed, one being the \emph{neighborhood grid}.
	
	In this paper, we consider combinatorial aspects of this data structure.
	Particularly, we show that an assumption on uniqueness, made in previous works, is not actually satisfied.
	We extend the notions of the neighborhood grid to arbitrary grid sizes and dimensions and provide two alternative, correct versions of the proof that was broken by the dissatisfied assumption.
	
	Furthermore, we explore both the uniqueness of certain states of the data structure as well as when the number of these states is maximized.
	We provide a partial classification by using the hook-length formula for rectangular Young tableaux.
	Finally, we conjecture how to extend this to all 2-dimensional cases.
\end{abstract}

\section{Introduction}
\label{sec:Motivation}

One of the most fundamental tasks in data processing, modeling, or simulations is finding answers to nearest neighbor queries.
These queries are to identify those points from a given data set that are closest or close enough to a certain query point.
Application examples include particle-, cell-, or crowd-simulations, where the movement of the particles, cells, or humans in the next time step is determined by the current placement and movement of their respective nearest neighbors.
As these queries have to be answered with respect to some embedding space, so-called \emph{spatial data structures} are designed to answer them.
The main metric to be optimized for these answers is time, as quick answers allow more simulation steps.

The \emph{Monotonic Logical Grid}~(MLG)\footnote{Some authors refer to it as the \emph{Monotonic Lagrangian Grid}.} is one such data structure that was developed specifically for large-scale simulations of molecular and fluid dynamics systems~\cite{boris1986vectorized}.
As opposed to, e.g., \emph{k-d trees}~\cite{skrodzki2019k-dtree}, the MLG does not answer neighborhood queries exactly.
That is, the MLG only provides approximate answers of neighbors that are somewhat close to the query point while not necessarily being the closest.
For the MLG, no guarantee can be given on the quality of the provided neighborhoods.
However, it can be implemented easily on highly-parallelized hardware and therefore computes its answers faster than, e.g., k-d trees, making it a suitable alternative for applications where speed trumps exactness.
See the article of Uhlmann~\cite{uhlmann2020monotonic} for an overview of the history of the MLG.

More recently, the MLG was rediscovered by researchers who were investigating crowd simulations~\cite{joselli2009neighborhood} and fluid animations~\cite{joselli2015neighborhood}.
They coined the data structure \emph{Neighborhood Grid}, a terminology which we will for the remainder of this paper. 
Their interest was sparked by the need for fast answers in their simulations.
Hence, a highly-parallelized \emph{CUDA} implementation on a graphics card that provided fast, approximate answers was exactly what was necessary for their applications.

To illustrate the neighborhood grid and how it is used in an application scenario, we will give a brief example.
A more formalized description is given in Section~\ref{sec:DefiningPointSetsAndStableStates}.
\begin{example}\label{exa:ExamplePointSet}
    Consider the set of 16 two-dimensional points in Euclidean space
    \begin{align*}
        A(2.0,3.1), B(2.5,3.4), C(2.9,3.2), D(2.2,4.7), E(1.4,4.4), F(0.5,5.0), G(1.0,3.6), H(4.1,4.8),\\ I(0.4,0.7), J(1.1,2.0), K(4.2,3.8), L(1.7,1.9), M(3.5,1.8), N(2.1,1.5), O(1.9,0.6), P(3.3,0.8).
    \end{align*}
    This point set is shown in the left of Figure~\ref{fig:ExamplePointSet}.
    
    To estimate neighborhoods of points in this point set, the neighborhood grid for these points is built by moving them into a~$4\times4$ grid, with each grid cell holding exactly one point.
    Furthermore, the points are arranged in a way such that the $x$-coordinate in the rows increases from left to right and the $y$-coordinates in the columns increases from bottom to top.
    A neighborhood grid like this is said to be in a \emph{stable state} (see Definition~\ref{def:StableState}).
    In Theorem~\ref{the:StableStateExistence}, we will establish that such a stable state always exists.
    See the right image of Figure~\ref{fig:ExamplePointSet} for an illustration of a stable state for the point set given above.

    To estimate the neighborhood of a query point from the point set (say, point~$A$), the one-ring around the query point in the neighborhood grid is considered.
    As an estimate, the nearest neighbor from this one-ring is output as the nearest neighbor of the query point.
    In the right of Figure~\ref{fig:ExamplePointSet}, the one-ring around query point~$A$ is shown in gray.
    The closest point to~$A$ in this one-ring is point $C$.
    As can be seen in the left of Figure~\ref{fig:ExamplePointSet}, the actual closest point to $A$ within the set is point~$B$.
    However, point~$C$ is a reasonable estimate for a nearest neighbor to~$A$, indicated by the dashed line.
\end{example}

\begin{figure}
    \centering
    \begin{tikzpicture}
        \draw[dashed, black!50] (2.0, 3.1) -- (2.9, 3.2);
        \draw[->] (-0.2,0) -- (5.0,0);
        \draw[->] (0,-0.2) -- (0,5.0);
        \node[] (x) at (2.5,-0.3) {$x$};
        \node[] (x) at (-0.3,2.5) {$y$};
        \draw[black!50] (2.0,3.1) circle (0.5831);
        \node[circle,fill=black,inner sep=2pt,label=left :{A}] (A) at (2.0, 3.1) {};
        \node[circle,fill=black,inner sep=2pt,label=above:{B}] (B) at (2.5, 3.4) {};
        \node[circle,fill=black,inner sep=2pt,label=right:{C}] (C) at (2.9, 3.2) {};
        \node[circle,fill=black,inner sep=2pt,label=right:{D}] (D) at (2.2, 4.7) {};
        \node[circle,fill=black,inner sep=2pt,label=above:{E}] (E) at (1.4, 4.4) {};
        \node[circle,fill=black,inner sep=2pt,label=below:{F}] (F) at (0.5, 5.0) {};
        \node[circle,fill=black,inner sep=2pt,label=left :{G}] (G) at (1.0, 3.6) {};
        \node[circle,fill=black,inner sep=2pt,label=right:{H}] (H) at (4.1, 4.8) {};
        \node[circle,fill=black,inner sep=2pt,label=below:{I}] (I) at (0.4, 0.7) {};
        \node[circle,fill=black,inner sep=2pt,label=left :{J}] (J) at (1.1, 2.0) {};
        \node[circle,fill=black,inner sep=2pt,label=right:{K}] (K) at (4.2, 3.8) {};
        \node[circle,fill=black,inner sep=2pt,label=above:{L}] (L) at (1.7, 1.9) {};
        \node[circle,fill=black,inner sep=2pt,label=right:{M}] (M) at (3.5, 1.8) {};
        \node[circle,fill=black,inner sep=2pt,label=right:{N}] (N) at (2.1, 1.5) {};
        \node[circle,fill=black,inner sep=2pt,label=below:{O}] (O) at (1.9, 0.6) {};
        \node[circle,fill=black,inner sep=2pt,label=below:{P}] (P) at (3.3, 0.8) {};
    \end{tikzpicture}
    \hspace{1.5cm}
    \begin{tikzpicture}[scale=0.71]
        \draw[fill=black!10] (0,2) rectangle (6,8);
        \draw[fill=white] (2,4) rectangle (4,6);
        \draw[step=2.0] (0,0) grid (8,8);
        \draw (3,5) node[scale=0.9] {A(2.0,3.1)}; 
        \draw (7,3) node[scale=0.9] {B(2.5,3.4)}; 
        \draw (5,5) node[scale=0.9] {C(2.9,3.2)}; 
        \draw (5,7) node[scale=0.9] {D(2.2,4.7)}; 
        \draw (3,7) node[scale=0.9] {E(1.4,4.4)}; 
        \draw (1,7) node[scale=0.9] {F(0.5,5.0)}; 
        \draw (1,5) node[scale=0.9] {G(1.0,3.6)}; 
        \draw (7,7) node[scale=0.9] {H(4.1,4.8)}; 
        \draw (1,1) node[scale=0.9] {I(0.4,0.7)}; 
        \draw (1,3) node[scale=0.9] {J(1.1,2.0)}; 
        \draw (7,5) node[scale=0.9] {K(4.2,3.8)}; 
        \draw (3,3) node[scale=0.9] {L(1.7,1.9)}; 
        \draw (7,1) node[scale=0.9] {M(3.5,1.8)}; 
        \draw (5,3) node[scale=0.9] {N(2.1,1.5)}; 
        \draw (3,1) node[scale=0.9] {O(1.9,0.6)}; 
        \draw (5,1) node[scale=0.9] {P(3.3,0.8)}; 
\end{tikzpicture}
    \caption{
        Sixteen points $A$ to $P$ in two-dimensional Euclidean space, note how the closest point to query point~$A$ is point~$B$. 
        The points are moved into a~$4\times 4$ grid in such a way that the rows of the grid are increasing in $x$-value from left to right while the columns of the grid are increasing in $y$-value from bottom to top, which forms a \emph{stable state} (see Definition~\ref{def:StableState}.
        Now, the closest point in the one-ring around a query point (a gray one-ring around query point~$A$ is shown) is taken as approximate nearest neighbor to the query point.
        Here, the closest point in the one-ring is~$C$, which is not the actual nearest neighbor~$B$, but a good approximation.
    }
    \label{fig:ExamplePointSet}
\end{figure}

Note that if the points are distinct when reducing to each dimension, the actual geometric coordinates of the points are not necessary to provide the configuration in the grid.
That is because we can order the points along a specific coordinate axis and replace this coordinate entry by the point's index in this order (which is unique if and only if the coordinates along this axis are distinct), starting from~$1$.
This operation reduces all coordinate entries to natural numbers, but does not change any relation of the points with respect to the grid representation.
This observation allows us to reduce a set of~$N$ geometric points to permutations of the first~$[N]$ natural numbers in each coordinate.
Hence, we can study all possibly arising neighborhood grids by reducing to their combinatorics.

We are not the first to start from a data structure and investigate it via its underlying combinatorics.
In fact, this approach has been applied almost throughout the entire history of data structure research, see, e.g.,~\cite{demillo1978small,woo1985combinatorial,wolter1991combinatorial,malashina2021combinatorial}.
It is a fruitful endeavour to consider computer science from the special standpoint of combinatorics, to the point that whole research projects can solely be devoted to this\footnote{Confer to the results of European Project 678765, \url{https://cordis.europa.eu/project/id/678765}.}.
Thereby, we are motivated to apply the combinatorial reduction outlined above to the neighborhood grid data structure. 
We will illustrate this reduction procedure on the previous example.

\begin{example}
    Consider the point set as given in Example~\ref{exa:ExamplePointSet}.
    Ordering the points by their first coordinate yields the order
    \begin{align*}
        (I, F, G, J, E, L, O, A, N, D, B, C, P, M, H, K),
    \end{align*}
    while sorting along the second coordinate yields the order
    \begin{align*}
        (O, I, P, N, M, L, J, A, C, B, G, K, E, D, H, F).
    \end{align*}
    Replacing the actual coordinates by their indices in these sortings then gives the following alternate point set
    \begin{align*}
        A(8,8),\ B(11,10),\ C(12,9),\ D(10,4),\ E(5,13),\ F(2,16),\ G(3,11),\ H(15,15),\ I(1,2),\ J(4,7),\\
        K(16,12),\ L(6,6),\ M(14,5),\ N(9,4),\ O(7,1),\ P(13,3).
    \end{align*}
    As shown in Figure~\ref{fig:ExamplePointSetReplaced}, replacing the original points in the neighborhood grid by this alternate point set keeps the neighborhood grid in a stable state.
\end{example}

\begin{figure}
    \centering
    \begin{tikzpicture}[scale=0.71]
        \draw[step=2.0] (0,0) grid (8,8);
        \draw (3,5) node[scale=0.9] {A(2.0,3.1)}; 
        \draw (7,3) node[scale=0.9] {B(2.5,3.4)}; 
        \draw (5,5) node[scale=0.9] {C(2.9,3.2)}; 
        \draw (5,7) node[scale=0.9] {D(2.2,4.7)}; 
        \draw (3,7) node[scale=0.9] {E(1.4,4.4)}; 
        \draw (1,7) node[scale=0.9] {F(0.5,5.0)}; 
        \draw (1,5) node[scale=0.9] {G(1.0,3.6)}; 
        \draw (7,7) node[scale=0.9] {H(4.1,4.8)}; 
        \draw (1,1) node[scale=0.9] {I(0.4,0.7)}; 
        \draw (1,3) node[scale=0.9] {J(1.1,2.0)}; 
        \draw (7,5) node[scale=0.9] {K(4.2,3.8)}; 
        \draw (3,3) node[scale=0.9] {L(1.7,1.9)}; 
        \draw (7,1) node[scale=0.9] {M(3.5,1.8)}; 
        \draw (5,3) node[scale=0.9] {M(2.1,1.5)}; 
        \draw (3,1) node[scale=0.9] {O(1.9,0.6)}; 
        \draw (5,1) node[scale=0.9] {P(3.3,0.8)}; 
    \end{tikzpicture}
    \hspace{1.5cm}
    \begin{tikzpicture}[scale=0.71]
        \draw[step=2.0] (0,0) grid (8,8);
        \draw (3,5) node[scale=0.9] {A(8,8)}; 
        \draw (7,3) node[scale=0.9] {B(11,10)}; 
        \draw (5,5) node[scale=0.9] {C(12,9)}; 
        \draw (5,7) node[scale=0.9] {D(10,14)}; 
        \draw (3,7) node[scale=0.9] {E(5,13)}; 
        \draw (1,7) node[scale=0.9] {F(2,16)}; 
        \draw (1,5) node[scale=0.9] {G(3,11)}; 
        \draw (7,7) node[scale=0.9] {H(15,15)}; 
        \draw (1,1) node[scale=0.9] {I(1,2)}; 
        \draw (1,3) node[scale=0.9] {J(4,7)}; 
        \draw (7,5) node[scale=0.9] {K(16,12)}; 
        \draw (3,3) node[scale=0.9] {L(6,6)}; 
        \draw (7,1) node[scale=0.9] {M(14,5)}; 
        \draw (5,3) node[scale=0.9] {N(9,4)}; 
        \draw (3,1) node[scale=0.9] {O(7,1)}; 
        \draw (5,1) node[scale=0.9] {P(13,3)}; 
    \end{tikzpicture}
    \caption{
        Replacing the actual point coordinate entries by their indices of a sorting in this coordinate keeps the property of a stable state in the neighborhood grid.
    }
    \label{fig:ExamplePointSetReplaced}
\end{figure}

Regarding the existince of stable states, Malheiros and Walter provide an algorithm to build these for arbitrary point sets~\cite{malheiros2015simple}.
However, their proof is given solely for two-dimensional point sets.
Furthermore, it relies on the implicit assumption that each point set has a unique stable state.
Building on the combinatorical viewpoint on the neighborhood grid and on this previous work, in this paper, we will:
\begin{itemize}
    \item \ldots show that any $d$-dimensional point set has a stable state, even for grids with varying side lengths. These stable states can be constructed in polynomial time (Theorem~\ref{the:StableStateExistence}).
    \item \ldots prove a lower bound on the run time of this construction (Theorem~\ref{the:lowerBound}) and prove time optimality of the construction algorithm for the case where all side lengths of the grid are equal (Corollary~\ref{cor:timeOptimality}).
    \item \ldots disprove the implicit assumption of Malheiros and Walter by showing that all point sets placed in grids of size~$4\times 4$ or larger admit to at least two stable states (Propositions~\ref{pro:U(4,4)=0} and~\ref{pro:PropagatingNon-Uniqueness}).
    \item \ldots bound the number of stable states that a point set admits to (Proposition~\ref{prop:youngnum}).
    \item \ldots conjecture a classification for the point sets with maximally many stable states for all 2-dimensional grids (Conjecture~\ref{conj:fibon}).
\end{itemize}
Preliminary results of this paper were published in a PhD thesis~\cite{skrodzki2019neighborhood} and presented at the EuroCG2018 conference~\cite{skrodzki2018combinatorial}.


\section{Defining Point Sets and Stable States}
\label{sec:DefiningPointSetsAndStableStates}

For natural numbers~$d$ and~$N$, let $\bm{\sigma} = (\sigma_1,\dots,\sigma_d)$ be a collection of permutations of the set~$[N]$ (i.e., each~${\sigma_i \in \mathfrak S_N}$). 
For each~${i \in [N]}$, we write ${\bm{\sigma}(i) = (\sigma_1(i),\sigma_2(i),\dots,\sigma_d(i))}$. 
We call the collection $P_{\bm\sigma} \coloneqq \bigcup_{i=1}^N \bm{\sigma}(i)$ a \emph{point set}, which we think of as the relative positions of a collection of points in $d$-space. 
We will often set $\sigma_1 = \identityperm$, where $\identityperm$ is the identity permutation, so that every point set is associated with a unique choice of permutations. 

Let~${(n_1,\dots,n_d) \in \mathbb N^d}$ such that~${\prod_{i=1}^d n_i = N}$. 
Consider the tuple~${(a_1,\dots,a_d)}$ such that~${1 \le a_i \le n_i}$ for all~${i \in [1,d]}$.
We say that these tuples make up the \emph{positions} on the $(n_1\times \ldots \times n_d)$ \emph{grid}.
For~${(k_1,\ldots,k_d)\in\mathbb{N}^d}$ such that~${k_i\leq n_i}$ for all~${i\in[d]}$ and for some position~${\mathbf{c} = (c_1,\ldots,c_d)}$ such that~${c_i+k_i \leq n_i}$, we call the set of positions
\begin{align*}
    {\{(b_1,\ldots,b_d)\mid c_i\leq b_i\leq c_i+k_i, i\in[d]\}}
\end{align*}
a~${k_1\times\ldots\times k_d}$ \emph{connected subgrid anchored at $\mathbf{c}$}.

Let $\pi$ be a bijective map $\pi: [n_1]\times\ldots\times[n_d] \to N$ which we call a \emph{placement map}.
Given collection~${\bm{\sigma} = (\sigma_1,\dots,\sigma_d)}$, let $M_\pi(P_{\bm\sigma})$ be the~${(n_1\times\ldots \times n_d)}$ grid where the \emph{position}~${\bm{a} = (a_1,\dots,a_d)}$ is labeled with the \emph{point}
\begin{align*}
    \bm{\sigma}(\pi(\bm{a})):= \left (\sigma_1(\pi(a_1,\dots,a_d)),\sigma_2(\pi(a_1,\dots,a_d)), \dots, \sigma_d(\pi(a_1,\dots,a_d))\right).
\end{align*}
We call~$M_\pi(P_{\bm\sigma})$ a \emph{state}. 
See Figure~\ref{fig:3states} for three examples of states for a given point set.

\begin{figure}
\centering
    \begin{tikzpicture}[scale=0.6]
        \draw[step=2.0] (0,0) grid (6,6);
        \draw (1,1) node[scale=0.9] {(1,3)}; 
        \draw (3,1) node[scale=0.9] {(2,2)}; 
        \draw (5,1) node[scale=0.9] {(8,1)}; 
        \draw (1,3) node[scale=0.9] {(3,4)}; 
        \draw (3,3) node[scale=0.9] {(4,6)}; 
        \draw (5,3) node[scale=0.9] {(5,5)}; 
        \draw (1,5) node[scale=0.9] {(6,8)}; 
        \draw (3,5) node[scale=0.9] {(7,7)}; 
        \draw (5,5) node[scale=0.9] {(9,9)}; 
        \draw (3,7) node[scale=0.9] {Stable State};
    \end{tikzpicture}
    \hspace{.5 cm}
    \begin{tikzpicture}[scale=0.6]
            \draw[step=2.0] (0,0) grid (6,6);
            \draw (1,1) node[scale=0.9] {(2,2)}; 
            \draw (3,1) node[scale=0.9] {(5,5)}; 
            \draw (5,1) node[scale=0.9] {(8,1)}; 
            \draw (1,3) node[scale=0.9] {(1,3)}; 
            \draw (3,3) node[scale=0.9] {(4,6)}; 
            \draw (5,3) node[scale=0.9] {(6,8)}; 
            \draw (1,5) node[scale=0.9] {(3,4)}; 
            \draw (3,5) node[scale=0.9] {(7,7)}; 
            \draw (5,5) node[scale=0.9] {(9,9)}; 
            \draw (3,7) node {Stable State};
    \end{tikzpicture}
    \hspace{.5 cm}
    \begin{tikzpicture}[scale=0.6]
        \draw[step=2.0] (0,0) grid (6,6);
        \draw (1,1) node[scale=0.9] {(2,2)}; 
        \draw (3,1) node[scale=0.9] {(5,5)}; 
        \draw (5,1) node[scale=0.9] {(4,6)}; 
        \draw (1,3) node[scale=0.9] {(1,3)}; 
        \draw (3,3) node[scale=0.9] {(7,7)}; 
        \draw (5,3) node[scale=0.9] {(8,1)}; 
        \draw (1,5) node[scale=0.9] {(3,4)}; 
        \draw (3,5) node[scale=0.9] {(6,8)}; 
        \draw (5,5) node[scale=0.9] {(9,9)}; 
        \draw (3,7) node {Non-Stable State};
    \end{tikzpicture}
    \caption{
        For~$d=2$ and~$N=9$, consider the collection of two permutations~${\bm{\sigma} = (\identityperm, \langle 3,2,4,6,5,8,7,1,9\rangle)}$. 
        Above are three states for~$P_{\bm\sigma}$ on the $3\times 3$ grid. 
        The first two states are stable while the third is not as, e.g., the lowest row is not increasing in the first coordinate.
        As~${\sigma_1=\identityperm}$, we can easily read off the placement maps from the states, as the first coordinate defines the placement.
        For instance, in the leftmost state, we have: $\pi(1,1)=1$, $\pi(2,1)=2$, $\pi(2,2)=4$, etc.
    }
    \label{fig:3states}
\end{figure}

We are interested in finding a particular class of states called \emph{stable states}. 
Let~${\bm{a} = (a_1,\dots,a_d)}$ and~${\bm{b} = (b_1,\dots,b_d)}$ be two distinct grid positions in~$M_\pi(P_{\bm\sigma})$.
They address both the position in the grid underlying~$M_\pi(P_{\bm\sigma})$ as well as the labeling induced at said positions by~$\pi$ and~$\bm\sigma$.
We say that~$\bm{a}$ is \emph{compatible} with~$\bm{b}$ with respect to~$\pi$ if any of the following hold:
\begin{itemize}
    \item For some~$1 \le i < j \le d$, we have~$a_i \not= b_i$ and~$a_j \not= b_j$.
    \item For some~$1 \le i \le d$, we have~$a_i < b_i$ and~$\sigma_i(\pi(a_1,\dots,a_d)) < \sigma_i(\pi(b_1,\dots,b_d))$.
    \item For some~$1 \le i \le d$, we have~$a_i > b_i$ and~$\sigma_i(\pi(a_1,\dots,a_d)) > \sigma_i(\pi(b_1,\dots,b_d))$.
\end{itemize}

\begin{definition}\label{def:StableState}
A state~$M_\pi(P_{\bm\sigma})$ is a \emph{stable state} if every pair of grid positions are compatible. 
\end{definition}

For each~${i \in [d]}$, an \emph{$i$-strip} is a maximal collection of positions in~${[n_1]\times\ldots\times[n_d]}$ which differ only in their~$i^{\text{th}}$ entry. 
In particular, for any fixed 
\begin{align*}
    (a_1,\dots,a_{i-1},a_{i+1},\dots,a_d) \in [n_1] \times \ldots \times [n_{i-1}] \times [n_{i+1}] \times \ldots \times [n_d], 
\end{align*}    
we get an $i$-strip of the form $\{(a_1,\dots,a_{i-1},k,a_{i+1},\dots a_d) \mid 1 \le k \le n_i\}$. 
Notice that $M_\pi(P_{\bm\sigma})$ is stable if and only if for every $1\le i \le d$, the labeling along each $i$-strip is increasing in its $i^{\text{th}}$ coordinate. 
We will also call $1$-strips \emph{rows} and $2$-strips \emph{columns}.
Conversely, for each~${i \in [d]}$, an \emph{$i$-slice} is a maximal collection of positions in~${[n_1]\times\ldots\times[n_d]}$ which can differ everywhere, but in their~$i^{\text{th}}$ entry.

\begin{example}
\label{exa:stableStates}
If~$d=2$, then~$M_\pi(P_{\bm\sigma})$ is a matrix whose entries are ordered pairs of integers.
In this case,~$M_\pi(P_{\bm\sigma})$ is a \emph{stable state} if the rows are ordered by first coordinate and the columns are ordered by second coordinate. 
See Figure~\ref{fig:3states} for an example of stable states and a non-stable state. 
\end{example}

\begin{lemma}
\label{lem:stablewelldef}
    If~$\bm\sigma$ and~$\bm{\sigma}'$ are collections of permutations such that~${P_{\bm\sigma} = P_{\bm{\sigma}'}}$, then the stable states of~$\bm\sigma$ and~$\bm{\sigma}'$ coincide. 
\end{lemma}
\begin{proof}
    Let~${\bm{\sigma} = (\sigma_1,\dots, \sigma_d)}$ and~${\bm{\sigma}' = (\sigma'_1,\dots, \sigma'_d)}$.
    Then, the condition that ${P_{\bm\sigma} = P_{\bm{\sigma}'}}$ is equivalent to the condition that there is a~${\bm{\sigma}^* \in \mathfrak S_N}$ such that for every~${i \in [d]}$, we have~$\sigma'_i = \sigma_i\circ \sigma^*$. 
    It is immediate that~$M_\pi(P_{\bm{\sigma}'})$ corresponds to~$M_{\bm{\sigma}^*\circ\pi }(P_{\bm\sigma})$. The lemma follows. 
\end{proof}

\begin{definition}
$\mathcal S(P_{\bm\sigma})$ is the set of all stable states for a point set~$P_{\bm\sigma}$.
\end{definition}
Notice that this set is well defined by Lemma~\ref{lem:stablewelldef}. 

\section{A Polynomial-Time Building Algorithm}
\label{sec:PolynomialTimeBuildingAlgorithm}

\noindent Starting from the definitions as presented in the previous section, we will now investigate the following question: Given any point set~$P_{\bm\sigma}$ as specified above, does it always hold that~$|\mathcal{S}(P_{\bm\sigma})|>0$, i.e., is there always at least one stable state~$M_\pi(P_{\bm\sigma})$?
Malheiros and Walter have answered this questions for the case of~$d=2$ and~$n_1=n_2$, see~\cite[Sec.~III, pp.~181--182]{malheiros2015simple}.
They did so by introducing a sorting algorithm that for any~$P_{\bm\sigma}$ produces a placement~$\pi$ such that~$M_\pi(P_{\bm\sigma})$ is a stable state.
In the following theorem, we will generalize this statement and their algorithm to the case of arbitrary grid lengths~$n_i$ and arbitrary dimension~$d$.

\begin{theorem}[generalizing~\cite{malheiros2015simple}]
	\label{the:StableStateExistence}
	For every natural number $N=\prod_{i=1}^d n_i$, $(n_1,\ldots,n_d)\in\mathbb{N}^d$, $d\in\mathbb{N}$, and for every collection of permutations ${\bm{\sigma}=(\sigma_1,\ldots,\sigma_d)}$, $\sigma_i\in\mathfrak{S}_N$, there is a placement map~$\pi$ such that $M_\pi(P_{\bm\sigma})$ is a stable state.
	Furthermore, this map~$\pi$ can be constructed in~${\mathcal{O}\left(N\sum_{j=0}^{d-1}\log\left(\prod_{i=1}^{d-j}n_i\right)\right)}$.
\end{theorem}
\begin{proof}
    Start with the points~$p_1,\ldots,p_N$ of~$P_{\bm\sigma}$,~$N=\prod_{i=1}^d n_i$,~$d\geq1$ interpreted as a first sequence.
    Use it as input for Algorithm~\ref{alg:ExistenceStableState}.
    The algorithms recursively sorts the input along its highest dimension and splits it into equally large subsequences, stripping them of their highest dimension entry, to be processed in the next recursion step.
    From the returned, ordered subsequences, a grid~$M$ is assembled.
    
    \begin{algorithm}
    	\caption{Recursive computation of a stable placement~$\pi$} 
    	\label{alg:ExistenceStableState}
    	\begin{algorithmic}[1]
    	    \Procedure{RecursiveSplit}{Sequence~$P'$, lengths $(n_1,\ldots,n_\ell)$}
    	        \State Sort $P'$ by the entries in the $\ell$th dimension
    	        \State Denote this ordering by $q_1,\ldots,q_{N'}$
    	        \If{$\ell=1$}
    	            \State \Return the sequence $(q_1,\ldots,q_{N'})$
    	        \Else
    	            \State $M \gets$ an empty $n_1\times\ldots\times n_\ell$ grid
    	            \State $N'' \gets \prod_{i=1}^{\ell-1} n_i$
    	            \For{$i=1,\ldots,n_\ell$}
    	                \State $Q_i \gets \{q_{(i-1)N''+1}\ldots,q_{iN''}\}$ \Comment{Split the sorted sequence into $n_\ell$ subsequences $Q_i$ of length $N''$ each}
    	                \State $M(\ldots,i)\gets$\Call{RecursiveSplit}{$Q_i$, $(n_1,\ldots,n_{\ell-1})$} \Comment{Fill each of the $n_\ell$ $i$-slices of the grid $M$ with a recursive solution}
    	            \EndFor
    	            \State \Return $M$
    	        \EndIf
    	    \EndProcedure
    	\end{algorithmic}
    \end{algorithm}
    
    The resulting grid~$M$ is then a stable state~$M_\pi(P_{\bm\sigma})$ by induction: If~${d=1}$, the grid is just a sorted sequence of numbers, which is trivially stable.
    For~${d>1}$, each of the ${(d-1)}$-dimensional slices of the resulting grid~$M$ is stable by the induction hypothesis. 
    Furthermore, it is stable along all $d$-strips as the glued sub-sequences have been ordered along this dimension in the first step of the procedure.
    Thus, the grid represents indeed a stable state.

    The runtime of the algorithm consists of the respective sorting steps, multiplied by the recursion level.
    It can be computed as
    \begin{align*}
        N\log(N) + (n_d) \left(\prod_{i=1}^{d-1}n_i\right)\log\left(\prod_{i=1}^{d-1}n_i\right) + \left(\prod_{j=d-1}^d n_j\right)\left(\prod_{i=1}^{d-2} n_i\right)\log\left(\prod_{i=1}^{d-2} n_i\right) +\\ \ldots + \left(\prod_{j=2}^d n_j\right)(n_1)\log(n_1)\\
        = N\sum_{j=0}^{d-1}\log\left(\prod_{i=1}^{d-j}n_i\right)
    \end{align*}
\end{proof}

Note that the implicit decision in the first step of the algorithm, to sort along the last coordinate, is arbitrary.
The theorem holds for sorting along any dimension, when corresponding splitting and filling is applied.
We will explore this observation in Section~\ref{sec:UniquenessOfStableStates}.

Note further that in the specific case of $n_i=n\in\mathbb{N}$ for all $i$, we have 
\begin{align*}
      & n^0(n^d\log(n^d)) + n^1(n^{d-1}\log(n^{d-1}) + \ldots + n^{d-1}(n^1\log(n))\\
    = & d(n^d\log(n)) + (d-1)(n^d\log(n)) + \ldots (d-(d-1))(n^d\log(n))\\
    = & (n^d\log(n))\sum_{i=1}^d i\\
    = & (n^d\log(n))\frac{d(d+1)}{2}
\end{align*}
When treating~$d$ as a constant, the stable state for this arrangement can therefore be found in time~${\mathcal{O}(n^d\log(n))}$.

\begin{example}
For the case~${d=2}$ and~${n=3}$, an illustration of the procedure presented by this theorem is given in Figure~\ref{fig:IllustrateSortingAlgorithm}.
\end{example}

\begin{figure}
	\centering
	\begin{framed}
		\scalebox{0.9}{
		\begin{minipage}{1.\textwidth}
			\centering
                $\left\{
                    (4,7),\ (3,8),\ (9,2),\ (7,5),\ (2,3),\ (1,6),\ (8,9),\ (5,1),\ (6,4)
                \right\}$\\[0.25cm]
			$\downarrow$ Consider all points as one sequence.~$\downarrow$\\[0.25cm]
			\begin{tikzpicture}[scale=0.6]
			\draw[step=2.0] (0,0) grid (18,2);
			\draw (1,1) node[scale=0.9] {(4,7)}; 
			\draw (3,1) node[scale=0.9] {(3,8)}; 
			\draw (5,1) node[scale=0.9] {(9,2)}; 
			\draw (7,1) node[scale=0.9] {(7,5)}; 
			\draw (9,1) node[scale=0.9] {(2,3)}; 
			\draw (11,1) node[scale=0.9] {(1,6)}; 
			\draw (13,1) node[scale=0.9] {(8,9)}; 
			\draw (15,1) node[scale=0.9] {(5,1)}; 
			\draw (17,1) node[scale=0.9] {(6,4)}; 
			\end{tikzpicture}\\[0.25cm]
			$\downarrow$ Sort the sequence according to the first coordinate.~$\downarrow$\\[0.25cm]
			\begin{tikzpicture}[scale=0.6]
			\draw[step=2.0] (0,0) grid (18,2);
			\draw (1,1) node[scale=0.9] {(1,6)}; 
			\draw (3,1) node[scale=0.9] {(2,3)}; 
			\draw (5,1) node[scale=0.9] {(3,8)}; 
			\draw (7,1) node[scale=0.9] {(4,7)}; 
			\draw (9,1) node[scale=0.9] {(5,1)}; 
			\draw (11,1) node[scale=0.9] {(6,4)}; 
			\draw (13,1) node[scale=0.9] {(7,5)}; 
			\draw (15,1) node[scale=0.9] {(8,9)}; 
			\draw (17,1) node[scale=0.9] {(9,2)}; 
			\end{tikzpicture}\\[0.25cm]	
			$\downarrow$ Separate into subsequences of size~$n=3$ and sort them according to the second coordinate.~$\downarrow$\\[0.25cm]
			\begin{tikzpicture}[scale=0.6]
			\draw[step=2.0] (0,0) grid (2,6);
			\draw[step=2.0] (4,0) grid (6,6);
			\draw[step=2.0] (8,0) grid (10,6);
			\draw (1,1) node[scale=0.9] {(2,3)}; 
			\draw (1,3) node[scale=0.9] {(1,6)}; 
			\draw (1,5) node[scale=0.9] {(3,8)}; 
			\draw (5,1) node[scale=0.9] {(5,1)}; 
			\draw (5,3) node[scale=0.9] {(6,4)}; 
			\draw (5,5) node[scale=0.9] {(4,7)}; 
			\draw (9,1) node[scale=0.9] {(9,2)}; 
			\draw (9,3) node[scale=0.9] {(7,5)}; 
			\draw (9,5) node[scale=0.9] {(8,9)}; 
			\end{tikzpicture}
		\end{minipage}
	}
	\end{framed}
	\caption{
	    An illustration of Algorithm~\ref{alg:ExistenceStableState} used in Theorem~\ref{the:StableStateExistence}. 
	    The last step gives the columns of the final matrix which is then in a stable state.
	}
	\label{fig:IllustrateSortingAlgorithm}
\end{figure}

Theorem~\ref{the:StableStateExistence} imposes an upper bound on the runtime of any time-optimal, comparison-based algorithm that creates a stable state $M_\pi(P_{\bm\sigma})$ for a point set~$P_{\bm\sigma}$. 
A next question to tackle is: What is a lower bound? 
A trivial lower bound is given by the number~$N$ of points in the point set. 
However, to sharpen this lower bound, we need further results on the cardinality of the set~${\mathcal{S}(P_{\bm\sigma})}$, i.e., on the number of stable states for a given point set~$P_{\bm\sigma}$. 


\section{Counting Stable States and a Lower Bound}
\label{sec:CombinatorialResults}

In order to better understand the search space for any sorting algorithm that is to build a stable state~$M_\pi(P_{\bm\sigma})$, we start by counting the number of possible stable states through straightforward combinatorial means. 

\begin{proposition}
	\label{the:Countings}
	Given some~${n_i\in\mathbb{N}}$ with~${N=\prod_{i=1}^d n_i}$, there are:
	\begin{enumerate}
		\item~$\left(N!\right)^{d-1}$ point sets~$P_{\bm\sigma}$ and~${\left(N!\right)^d}$ different states~$M_\pi(P_{\bm\sigma})$,
		\item~$(N!)^d/{\prod_{i=1}^d {(n_i!)^{N/n_i}}}$ ways to fill the grid with a stable state,
		\item~$N!/(n_1!)^{N/n_1}$ placements that are stable with respect to the first coordinate and each such placement is stable for~$(N!)^{d-1}/\prod_{i=2}^d{(n_i!)^{N/n_i}}$ point sets~$P_{\bm\sigma}$.
	\end{enumerate}
	Furthermore,
	\begin{enumerate}
		\item[(4)]~$1/\prod_{i=1}^d (n_i!)^{N/n_i}$ of all fillings of the grid are stable.
		\item[(5)] The expected number of stable states for a point set chosen uniformly at random is $N!/\prod_{i=1}^d (n_i!)^{N/n_i}$.
	\end{enumerate}
\end{proposition}
\begin{proof}
	For the first statement, recall that each point set~$P_{\bm\sigma}$ is created by~$d$ permutations~${\sigma_i\in\mathfrak{S}}$, where~$\sigma_1$ is the identity permutation~$\identityperm$. 
	Thus, are~$N!$ choices for~${2\leq i \leq d}$, yielding a total of~$(N!)^{d-1}$ possible point sets.
	Furthermore, there are~$N!$ ways to pick the map~$\pi$ which places the point set in the grid.
	Therefore, we obtain~$(N!)^d$ possible states~$M_\pi(P_{\bm\sigma})$.

    For the second statement, recall that when a dimension~${i\in[d]}$ is fixed, the $i^{\text{th}}$ coordinate has to be ascending along all $i$-strips.
    For a fixed~${i\in[d]}$, the number of possible partitions, i.e., the number of ways to fill the~$i$-strips with coordinates, is given by the multinomial coefficient
    \begin{align*}
        \genfrac(){0pt}{}{N}{n_i,\ldots,n_i}
        =\genfrac(){0pt}{}{N}{n_i} \genfrac(){0pt}{}{N-n_i}{n_i,\ldots,n_i}
        =\ldots
        =N!/(n_i!)^{N/n_i},
    \end{align*}
    because once the coordinates for a specific~$i$-strip are chosen, they have to be ordered increasingly, for which there is only one way.
    However, the choices for how to fill the~$i$-strips have to be made for each of the~${i\in d}$ dimensions. Therefore, the total number of stable states is given by \[\prod_{i=1}^d \frac{N!}{(n_i!)^{N/n_i}}=(N!)^d/\prod_{i=1}^d (n_i!)^{N/n_i}.\]

        The third statement follows from an analogous argument to the second statement. 
	
	For the fourth statement, we divide the number of stable states (given by statement 2) by the total number of states (given by statement 1). 

        Finally, the fifth statement follows from dividing the number of point sets by the number of stable states. 
\end{proof}

For the special case of~${n_i=n_j}$ for all~${i\neq j}$, we obtain the following corollary:
\begin{corollary}
	\label{cor:Countings}
	Given some~${N,n\in\mathbb{N}}$ with~${N=n^d}$, there are:
	\begin{enumerate}
		\item~$\left(N!\right)^{d-1}$ point sets~$P_{\bm\sigma}$ and therefore~${\left(N!\right)^d}$ different states~$M_\pi(P_{\bm\sigma})$,
		\item~${\left(N!/(n!)^{n^{d-1}}\right)^d}$ ways to fill the grid with a stable state,
		\item$\frac{N!}{(n!)^{d-1}}$ placements $\pi$ that are stable with respect to the first coordinate and each such placement is stable for ~${\left(N!/(n!)^{n^{d-1}}\right)^{d-1}}$ point sets~$P_{\bm\sigma}$.
	\end{enumerate}
	Furthermore,
	\begin{enumerate}
		\item[4.]~$1/(n!)^{n^{d(d-1)}}$ of all fillings of the grid are stable
		\item[5.] and each point set $P_{\bm\sigma}$ has $N!/(n!)^{dn^{d-1}}$ stable states on average.
	\end{enumerate}
\end{corollary}

Having the necessary results at hand, we can now prove a lower bound on the building time of the neighborhood grid. 

\begin{theorem}
	\label{the:lowerBound}
	Any comparison-based algorithm has to perform at least ${\omega\left(\log\left(\prod_{i=2}^d (n_i!)^{N/n_i})\right)\right)}$ operations.
\end{theorem}
\begin{proof}
    Consider any comparison-based algorithm~$\mathcal{A}$ that constructs a placement~$\pi$ such that~$M_\pi(P_{\bm\sigma})$ is a stable state for some given point set~$P_{\bm\sigma}$ with respect to the first coordinate.
    By Proposition~\ref{the:Countings}, there are~${\tfrac{N!}{(n_1!)^{N/n_1}}}$ many of these, where the algorithm has to identify the correct one.
    Each subsequent query of~$\mathcal{A}$ can be considered as a node of a decision-tree, where the leaves correspond to placements of which some are stable for~$P_{\bm\sigma}$. 
    For an optimal algorithm, this tree is balanced and has depth~$\log(N!)$. 
    Recall the result of Proposition~\ref{the:Countings} that any such placement~$\pi$ is stable for~${\prod_{i=2}^d \frac{N!}{(n_i!)^{N/n_i}}}$ point sets.
    Thus, when building the tree, the algorithm~$\mathcal{A}$ cannot stop at a subtree with more than
    \begin{align*}
        \tfrac{N!}{(n_1!)^{N/n_1}} \prod_{i=2}^d \frac{N!}{(n_i!)^{N/n_i}} = \prod_{i=1}^d \frac{N!}{(n_i!)^{N/n_i}}
    \end{align*}
    leaves, as either the placement will not be stable with respect to the first coordinate or one of the leaves will surely not be stable under the currently considered placement.
    Therefore, the algorithm has to perform at least 
    \begin{align*}
          \log\left(N!\right) - \log\left(\prod_{i=1}^d \frac{N!}{(n_i!)^{N/n_i}}\right)\\
        = \log(N!) - \left(\log(N!) + \log\left(\prod_{i=1}^d \frac{1}{(n_i!)^{N/n_i}}\right)\right)\\
        = \log\left(\prod_{i=1}^d (n_i!)^{N/n_i}\right)
    \end{align*}
    many comparisons before it can possibly terminate.
\end{proof}

For non-uniform configurations with~$n_i\neq n_j$ for some~$i\neq j$, it remains unclear whether a faster algorithm than that from Theorem~\ref{the:StableStateExistence} can be found.
However, Theorem~\ref{the:lowerBound} yields the following corollary for $d$-dimensional point sets of uniform shape~${n_i=n}$ for all~${i\in[d]}$ and some~${n\in\mathbb{N}}$.

\begin{corollary}
    \label{cor:timeOptimality}
    In the case of~${n_i=n\in\mathbb{N}}$ for all~${i\in[d]}$, Algorithm~\ref{alg:ExistenceStableState} is a time-optimal building algorithm for the neighborhood grid among all comparison-based algorithms and takes~$\Theta(n^d\log(n))$ time to build a grid for a point set~$P_{\bm{\sigma}}$ of~${N=n^d}$ points.
\end{corollary}
\begin{proof}
    In the case of~${n_i=n}$, the tree has to be traversed to depth at least
    \begin{align*}
          \log\left(\prod_{i=1}^d (n!)^{n^{d-1}}\right) 
        = \log\left(\left((n!)^{n^{d-1}}\right)^{d}\right) 
        = d\log\left((n!)^{n^{d-1}}\right)
        = dn^{d-1}\log(n!) = \mathcal{O}(n^d\log(n)).
    \end{align*}
    Therefore, each comparison-based algorithm building a stable state needs to perform at least~${\Omega(n^d\log(n))}$ operations. The observation that the property of a point set being restricted does not interfere with the run time completes the proof.
\end{proof}

Finally, we state the result in the form that it is historically most used in, cf.~\cite{malheiros2015simple}, namely the case of uniform grid side lengths in dimension~${d=2}$.
\begin{corollary}
    \label{cor:timeOptimalityD=2}
    In the case of~{$d=2$}, Algorithm~\ref{alg:ExistenceStableState} is a time-optimal building algorithm for the neighborhood grid among all comparison-based algorithm and takes~$\Theta(n^2\log(n))$ time to build a grid for a point set~$P_{\bm{\sigma}}$ of~${N=n^2}$ points.
\end{corollary}


\section{Point Sets with Unique Stable States}
\label{sec:UniquenessOfStableStates}

In the previous sections, we have constructively proven the existence of stable states for a given point set~$P_{\bm\sigma}$ on an~${n_1 \times \dots \times n_d}$ grid. 
However, in many cases, this stable state is not unique (as shown in Figure~\ref{fig:3states}). 
In this section, we will explore properties of point sets with unique stable states and partially classify the grids which allow for such point sets. 
We begin with two necessary conditions for a point set to have a unique stable state.

\begin{lemma}[First necessary condition of a unique stable state]
\label{lem:nc1}
    Consider a stable state~$M_\pi(P_{\bm\sigma})$ on a $d$-dimensional point set~$P_{\bm\sigma}$. 
    Suppose that~$M_\pi(P_{\bm\sigma})$ is the unique stable state of~$P_{\bm\sigma}$ and consider any~${i \in [d]}$. 
    For any grid positions~${\mathbf{a} = (a_1,\dots,a_d)}$ and~${\mathbf{b} = (b_1,\dots,b_d)}$ such that~${a_i < b_i}$, we must have~${\sigma_i(\pi(a)) < \sigma_i(\pi(b))}$.
\end{lemma}
\begin{proof}
    First, suppose that the condition is not satisfied for~${i=1}$. 
    Then, if we apply Theorem~\ref{the:StableStateExistence} to~$P_{\bm\sigma}$, we obtain a stable state which does satisfy the condition, because the algorithm begins by sorting by the first coordinate. 
    This contradicts the assumption that~$M_\pi(P_{\bm\sigma})$ is unique. 
    
    For general~$i$, we can consider an alternate version of the algorithm from Theorem~\ref{the:StableStateExistence} that switches the role of the first and~$i$-{th} coordinate. 
    Because the order of the coordinates is arbitrary, this will also produce a stable state, and our proof works equivalently to the~${i=1}$ case. 
\end{proof}

This provides a first necessary condition on the uniqueness of a stable state~$M_\pi(P_{\bm\sigma})$.
We can equivalently reformulate this condition on~$\bm\sigma$ in terms of a condition on~$\pi$ as follows.

\begin{corollary}
\label{cor:nc1}
    Given a point set~$P_{\bm\sigma}$ with~${N=\prod_{i=1}^d n_i}$ points and a stable state~$M_\pi(P_{\bm\sigma})$.
    In order for~$M_\pi(P_{\bm\sigma})$ to be the unique stable state of~$P_{\bm\sigma}$, it has to satisfy the following condition: For any given coordinate~${i\in [d]}$, the~$\prod_{j\in[d]\backslash i} n_j$ smallest $i$-coordinates have to be placed in the first positions of all~$i$-strips, the next~$\prod_{j\in[d]\backslash i} n_j$ smallest $i$-coordinates have to be placed in the second positions of all~$i$-strips, etc., until the~$\prod_{j\in[d]\backslash i} n_j$ largest $i$-coordinates are placed in the~$n_i$-th position of all~$i$-strips.
\end{corollary}

This first necessary condition provides candidates for point sets~$P_{\bm\sigma}$ that could have a unique stable placement.
We can count the number of point sets~$P_{\bm\sigma}$ that satisfy the condition of Lemma~\ref{lem:nc1} by considering the reformulation according to Corollary~\ref{cor:nc1}.
For each of the~${i\in[d]}$ coordinates, consider each of their~$n_i$ many~$i$-slices.
Each of these slices has~$\tfrac{N}{n_i}$ many points of which the~$i$-th coordinates can be permuted without violating the condition.
Therefore, there are
\begin{align*}
    \prod_{i\in [d]}\left(\left(\frac{N}{n_i}\right)!\right)^{n_i}
\end{align*}
many point sets that satisfy the condition of Lemma~\ref{lem:nc1}.
Throughout this section, we will have the case of~${d=2}$, ${n_1=n_2=4}$ as a running example.
In that case, there are~${16!\approx 2.1\cdot 10^{13}}$ point sets of which~${((4!)^4)^2\approx 1.1\cdot 10^{11}}$ satisfy the condition of Lemma~\ref{lem:nc1}.

To further simplify the challenge of checking which point sets have a unique stable state, we introduce a second necessary condition. This condition works recursively and allows results on smaller grids to be applied to larger grids. 

\begin{lemma}[Second necessary condition of a unique stable state]
\label{lem:nc2}
    Given a point set~$P_{\bm\sigma}$ with~${N=\prod_{i=1}^d n_i}$ points and a unique stable state~$M_\pi(P_{\bm\sigma})$. 
    Then, any~${k_1\times \ldots\times k_d}$ connected subgrid of~$M_\pi(P_{\bm\sigma})$ with~${k_i\in[d]}$ is in a unique stable state.
\end{lemma}
\begin{proof}
    Given some unique stable state~$M_\pi(P_{\bm\sigma})$. 
    Assume there exists some ${k_1\times\ldots\times k_d}$ connected subgrid~$\widetilde{M}$ of~$M$, anchored at~$\mathbf{c}=(c_1,\ldots,c_d)$, such that~$\widetilde{M}$ is not unique, but has a different stable state~$\overline{M}$. 
    
    By Corollary~\ref{cor:nc1}, we know that the~$i$-th values of the first~${c_i-1}$ elements in each $i$-strip are smaller than the $i$-th value of the $c_i$-th element.
    Furthermore, the $i$-th values of the last~${n_i-(c_i+k_i)}$ elements in each $i$-strip are larger than the $i$-th value of the~$c_i+k_i$ element.
    This fact remains true independent of any stable re-ordering of the~${k_1\times\ldots\times k_d}$ connected subgrid of~$M$, anchored at~$\mathbf{c}$.
    In particular, this remains true for the reordering induced by~$\overline{M}$.
    Therefore, replacing~$\widetilde{M}$ by~$\overline{M}$ in~$M$ gives another stable state, which violates the uniqueness of~$M$.
\end{proof}

For this statement, it is not as easy as it was for Lemma~\ref{lem:nc1} to give a count of those point sets that satisfy the condition.
This is mostly because of the recursive nature of the condition: we would not only have to take into account all connected subgrids of some~$M_\pi(P_{\bm\sigma})$, but also ensure that they interlock correctly.
However, we can turn to the case of~${d=2}$, ${n_1=n_2=4}$ again to provide a feeling of how Lemma~\ref{lem:nc2} reduces the number of point sets that can have a unique stable state.
Recall that there are approximately~${2.1\cdot10^{13}}$ point sets for this case, out of which approximately~${1.1\cdot10^{11}}$ satisfy Lemma~\ref{lem:nc1}.
We can put an upper bound on the number of point sets satisfying Lemma~\ref{lem:nc2} by considering all nine connected $2\times 2$~subgrids, for each of which there are~$12$ unique stable states.
Hence, there are at most~$12^9\approx 5.2\cdot10^9$ point sets that have a unique stable state.
While this brings the previous estimate down by two orders of magnitude, the actual number is far lower, as we will establish within the proof of Proposition~\ref{pro:U(4,4)=0}.

Note that the converses of Lemmas~\ref{lem:nc1} and~\ref{lem:nc2} do not hold: a stable state can satisfy both properties without being unique. 
For example, consider the point sets and their placements as shown in Figure~\ref{fig:NonEquivalenceXYBin}. 
The respective left stable states satisfy the conditions of Lemmas~\ref{lem:nc1} and~\ref{lem:nc2}, but the right state shows that these conditions do not imply uniqueness.

\begin{figure}
	\centering
	\begin{tikzpicture}[scale=0.6]
		\draw[step=2.0] (0,0) grid (4,4);
		\draw (1,3) node[scale=0.9] {(2,3)}; 
		\draw (3,3) node[scale=0.9] {(4,4)}; 
		\draw (1,1) node[scale=0.9] {(1,1)}; 
		\draw (3,1) node[scale=0.9] {(3,2)}; 
	\end{tikzpicture}
	\hspace{0.125cm}
	\begin{tikzpicture}[scale=0.6]
		\draw[step=2.0] (0,0) grid (4,4);
		\draw (1,3) node[scale=0.9] {(3,2)}; 
		\draw (3,3) node[scale=0.9] {(4,4)}; 
		\draw (1,1) node[scale=0.9] {(1,1)}; 
		\draw (3,1) node[scale=0.9] {(2,3)}; 
	\end{tikzpicture}
	\hspace{0.25cm}
	\begin{tikzpicture}[scale=0.6]
		\draw[step=2.0] (0,0) grid (6,6);
		\draw (1,1) node[scale=0.9] {(3,3)}; 
		\draw (1,3) node[scale=0.9] {(2,4)}; 
		\draw (1,5) node[scale=0.9] {(1,9)}; 
		\draw (3,1) node[scale=0.9] {(6,2)}; 
		\draw (3,3) node[scale=0.9] {(5,5)}; 
		\draw (3,5) node[scale=0.9] {(4,8)}; 
		\draw (5,1) node[scale=0.9] {(9,1)}; 
		\draw (5,3) node[scale=0.9] {(8,6)}; 
		\draw (5,5) node[scale=0.9] {(7,7)}; 
	\end{tikzpicture}
	\hspace{0.125cm}
	\begin{tikzpicture}[scale=0.6]
		\draw[step=2.0] (0,0) grid (6,6);
		\draw (1,1) node[scale=0.9] {(2,4)}; 
		\draw (1,3) node[scale=0.9] {(4,8)}; 
		\draw (1,5) node[scale=0.9] {(1,9)}; 
		\draw (3,1) node[scale=0.9] {(3,3)}; 
		\draw (3,3) node[scale=0.9] {(5,5)}; 
		\draw (3,5) node[scale=0.9] {(7,7)}; 
		\draw (5,1) node[scale=0.9] {(9,1)}; 
		\draw (5,3) node[scale=0.9] {(6,2)}; 
		\draw (5,5) node[scale=0.9] {(8,6)}; 
	\end{tikzpicture}
	\caption{
	    The left~$2\times 2$ and~$3\times 3$ stable states respectively satisfy the conditions of Lemmas~\ref{lem:nc1} and~\ref{lem:nc2}.
	    While the condition of Lemma~\ref{lem:nc2} is trivially satisfied for the~$2\times 2$ example, it has to be checked for the~$3\times 3$ example.
	    The right~$2\times 2$ and~$3\times 3$ stable states show that the conditions of Lemmas~\ref{lem:nc1} and~\ref{lem:nc2} are only necessary conditions, but do not imply uniqueness of the stable states.
	}
	\label{fig:NonEquivalenceXYBin}
\end{figure}

\begin{definition}
    Let $(n_1,\ldots,n_d)\in\mathbb{N}^d$, then we write $\mathcal{U}(n_1,\dots,n_d)$ for the number of point sets which have a unique stable state on the~${(n_1 \times \dots \times n_d)}$ grid. 
\end{definition}

To give some examples, it is easy to check that of the~$24$ different point sets on four points, it holds that~${\mathcal{U}(2,2)=12}$.
The other twelve point sets allow for two stable states each.
When placing nine points in a~$(3\times3)$ grid, checking the the candidates that satisfy Lemmas~\ref{lem:nc1} and~\ref{lem:nc2} reveals that~${\mathcal{U}(3,3)=966}$, out of~362,880 point sets.
Considering these numbers, one natural question to ask is whether there exist grids for which~${\mathcal{U}(n_1,\ldots,n_d)=0}$. 
In fact, we prove that there are infinitely many such grids in the next two propositions. 

\begin{proposition}
\label{pro:U(4,4)=0}
    $\mathcal{U}(4,4) = 0$
\end{proposition}
\begin{proof}
    Making use of the reductions provided by Lemmas~\ref{lem:nc1} and~\ref{lem:nc2}, we proved the statement using a brute force computational approach.
    For this, we assume that we have a complete list of unique stable states that fit all possible subgrids of the~${4\times 4}$ grid.
    These are 12, 106, and 966 for the~$2\times2$, $2\times 3$, and~$3\times3$ case, respectively.
    As these instances are small, they are easily computed and the time to do so is negligible within the overall time used to run the program.
    This brings the upper bound down to~$966\cdot106^2\cdot12\approx 1.3\cdot10^8$.
    
    Given this data for the subgrids, the program iteratively constructs all point sets that satisfy the condition of Lemma~\ref{lem:nc1}.
    For each point set, first, the condition of Lemma~\ref{lem:nc2} is checked and those that do not satisfy it are discarded.
    This leaves a mere~37,536 candidate point sets for the~$4\times4$ case.
    In order to optimally utilize Lemma~\ref{lem:nc2}, the construction of the point sets is done in a lazy way: Only part of the point set is constructed and if this already violates Lemma~\ref{lem:nc2}, all point sets containing this part can be neglected.
    This shows the importance of Lemma~\ref{lem:nc2}, as iterating over all approximately~$10^{8}$ point sets satisfying the condition of Lemma~\ref{lem:nc1} would cause a significantly longer computation time.
    Thus, by combining Lemmas~\ref{lem:nc1} and~\ref{lem:nc2}, we reduce the order of magnitude of point sets to check from the trivial bound of~$10^{11}$ to~$10^4$ for our example.
    
    Finally, for the remaining candidate point sets, all~$16!$ different placements of the point set in the~$(4\times 4)$ grid are checked.
    Again, this is only done until two stable placements are found, which contradicts uniqueness for this point set and the iteration over all other placements can be halted.
    
    After running this program, we have found two stable states for all point sets that satisfy the conditions of Lemma~\ref{lem:nc2}.
    Hence, we can conclude that there is not a single point set with a unique placement in the~$(4\times 4)$ grid.
\end{proof}

Hence, the~$(4\times4)$ grid does not allow for unique stable states.
We can extend this to other grids by the following proposition.

\begin{proposition}
\label{pro:PropagatingNon-Uniqueness}
    For two natural numbers~${d \geq d'}$, let $(n_1,\dots,n_d)\in\mathbb{N}^{d}$, $(n'_1,\dots,n'_{d'})\in\mathbb{N}^{d'}$ such that for all~${i \in [d']}$, we have~${n_i \geq n'_i}$.
    If~${\mathcal{U}(n'_1,\ldots,n'_{d'}) = 0}$, then~${\mathcal{U}(n_1,\ldots,n_d) = 0}$.
\end{proposition}
\begin{proof}
    Assume there is a stable state on an~$(n_1\times\ldots\times n_d)$ grid that satisfies the condition of Lemma~\ref{lem:nc1}.
    As~${n_i \geq n'_i}$ for all~$i\in[d']$, there exists an~${n'_1\times\ldots\times n'_{d'}\times 1\times \ldots\times 1}$ connected subgrid (with~${d-d'}$ trailing 1s) which has no unique stable state.
    This violates the condition of Lemma~\ref{lem:nc2}, which implies that the~$(n_1\times\ldots\times n_d)$ grid also does not have any stable states.
\end{proof}

By Proposition~\ref{pro:PropagatingNon-Uniqueness}, a grid that has any subgrid without unique stable states does not have unique stable states itself.
Combining this with Proposition~\ref{pro:U(4,4)=0}, we arrive at the following corollary.

\begin{corollary}
\label{cor:U(n_1,...,n_d)=0}
    Let~$(n_1,\ldots,n_d)\in\mathbb{N}^d$ that for some~${i \neq j}$, we have~${n_i,n_j \geq 4}$. 
    Then, we must have ${\mathcal{U}(n_1,\ldots,n_d) = 0}$.
\end{corollary}

This proves that there are infinitely many grids that do not admit a single unique stable state.
Going back to the paper of~\cite{malheiros2015simple}, as stated above, they provided a proof for an upper bound on the time needed to build a stable state, see the discussion before Theorem~\ref{the:lowerBound}.
In their paper, they continue to say that ``the problem of sorting $s$ unrelated lists of $s$ real values has~${\mathcal{O}(n log n)}$ as its established lower bound''~\cite[p.~182]{malheiros2015simple}, where they assume~$d=2$; their~$s$ would be $n_1,n_2$ and their~$n$ is~$N$ is our notation.
This serves as their argument to establish a lower bound that coincides with the upper one.

Note that this argument assumes that in the two-dimensional case, there is a unique separation of the points into the rows (or columns) of the grid.
Given this assumption, the rows (or columns) have a respective unique sorting.
However, our Corollary~\ref{cor:U(n_1,...,n_d)=0} shows that each point set has at least two stable states.
Each of these stable states provides a different partition of the points into the rows (or columns).
Thereby, the assumption of unique sets of rows (or columns) is incorrect.
Thus, the argument for the existence of a lower bound as provided by~\cite{malheiros2015simple} does not hold.

Nevertheless, the proof as given by us for Corollary~\ref{cor:timeOptimality} does not depend on any assumption of uniqueness.
Thus, this proof of time-optimality -- for the case of~$n_i=n$ -- repairs the statement of Malheiros and Walter and even extends it to arbitrary grid shapes.

To continue, note that Corollary~\ref{cor:U(n_1,...,n_d)=0} proves that there are infinitely many point sets without a unique stable placement.
Now, we turn to the converse and ask: How many grids do allow for unique stable states?
In fact, we also give an infinite family of grids such that~${\mathcal{U}(n_1,\ldots,n_d) > 0}$.

\begin{proposition}
\label{prop:2unique}
    If~${n_i = 2}$ for all~${i<d}$ (with~${n_d\in\mathbb{N}}$ arbitrary), then~${\mathcal{U}(n_1,\ldots,n_d) > 0}$. 
\end{proposition}
\begin{proof}
    We can explicitly construct a point set with a unique stable state. 
    Recall that 
    \begin{align*}
        N = \prod_{i=1}^d n_i = 2^{d-1}n_d.
    \end{align*}
    Let~${\sigma_1 = \identityperm}$ and for~${i>1}$, define:
    \begin{align*}
        \sigma_i(j) = 
        \begin{cases}
            \min k\ge 1 \text{ such that } k \not\in \{\sigma_i(1),\dots,\sigma_i(j-1)\} & \text{ if $\lfloor j/{2^{i-1}}\rfloor$ is even},\\
            \max k\le N \text{ such that } k \not\in \{\sigma_i(1),\dots,\sigma_i(j-1)\} & \text{ if $\lfloor j/2^{i-1}\rfloor$ is odd}.
        \end{cases}
    \end{align*}
    We then place the points in the grid sequentially from~$p_1$ up to~$p_N$. 
    Notice that for each~${k \in [N]}$, once we have placed the points~${p_1,\ldots,p_k}$ into the grid, there is only one place where $p_{k+1}$ can go without forcing the state to be non-stable.
    
    This can be seen by inductively considering the $d$-slices of the grid to be filled.
    Arriving at an empty $d$-slice, there are~$2^{d-1}$ points to be placed.
    The next~$2^{d-1}$ points given by~$\bm\sigma$ have the property that each of their coordinates is either minimal or maximal with respect to the currently available numbers.
    Furthermore, among these next~$2^{d-1}$ points, each of the~$2^{d-1}$ many choices of minimum and maximum for each coordinate occurs exactly once.
    This places the points uniquely within the currently considered $d$-slice.
    
    Furthermore, because the minimum or maximum is chosen, the first~$2^{d-1}$ points of~$\bm\sigma$ have to go into the first $d$-slice, the next~$2^{d-1}$ points have to go into the second slice, and so on, as otherwise there would be a pair of points along an $i$-strip, $1\leq i\leq d-1$ that is not stably placed.
    Hence, there is only one stable placement for this point set and ${\mathcal{U}(2,\ldots,2,n_d)\geq1}$.
\end{proof}

\begin{example}
    For~${1 \leq i \leq 5}$, the construction of Proposition~\ref{prop:2unique} gives the following permutations:
    \begin{align*} 
        \sigma_1 &= (1,2,3,4,5,6,7,8,9,10,\dots,N),\\
        \sigma_2 &= (1,N,2,N-1,3,N-2,4,N-3,5,N-4,\dots),\\
        \sigma_3 &= (1,2,N,N-1,3,4,N-2,N-3,5,6,\dots),\\
        \sigma_4 &= (1,2,3,4,N,N-1,N-2,N-3,5,6,\dots),\\
        \sigma_5 &= (1,2,3,4,5,6,7,8,N,N-1,\dots).
    \end{align*}
    See Figure~\ref{fig:unique2byn} for the case where $d=2$, $n_1 = 2$, and $n_2 = 5$.
    Note that all points with odd first coordinates have to go into the bottom row because of their small second coordinate.
    Similarly, all points with even first coordinates have to go into the upper row because of their large second coordinate.
    Thus, there is a unique way to place each point while maintaining stability. 
\end{example}

\begin{figure}
		\centering
		\begin{tikzpicture}[scale=0.6]
		\draw[step=2.0] (0,0) grid (10,4);
		\draw (1,1) node[scale=0.9] {(1,1)}; 
		\draw (1,3) node[scale=0.9] {(2,10)}; 
		\draw (3,1) node[scale=0.9] {(3,2)}; 
		\draw (3,3) node[scale=0.9] {(4,9)}; 
		\draw (5,1) node[scale=0.9] {(5,3)}; 
		\draw (5,3) node[scale=0.9] {(6,8)}; 
		\draw (7,1) node[scale=0.9] {(7,4)}; 
		\draw (7,3) node[scale=0.9] {(8,7)}; 
		\draw (9,1) node[scale=0.9] {(9,5)}; 
		\draw (9,3) node[scale=0.9] {(10,6)}; 

		\end{tikzpicture}
		\caption{
		    For $d = 2$, $n_1 = 2$, and $n_2 = 5$, Proposition~\ref{prop:2unique} says that there exists a point set with a unique sable state. 
		    Above is the the unique stable state for the point set $P_{(\sigma_1,\sigma_2)}$ as constructed in the proof of the proposition.
		}
		\label{fig:unique2byn}
\end{figure}

For the remainder of this section, we will restrict to the case of~${d=2}$. From Proposition~\ref{pro:U(4,4)=0} and Corollary~\ref{cor:U(n_1,...,n_d)=0}, we can deduce that~${\mathcal{U}(n_1,n_2)=0}$ for all~${n_1,n_2\geq 4}$. Thus, it suffices to restrict to the case where $n_1 \in [3]$. 

When~${n_1=1}$, the number of unique stable states is just the number of point sets, as these are a simple sequence and are therefore uniquely ordered. Hence, we have~${\mathcal{U}(1,n_2)=n_2!}$.
When~${n_1=2}$, Proposition~\ref{prop:2unique} implies that ${\mathcal{U}(2,n_2)} > 0$. 
When~${n_1=3}$, we conjecture that ${\mathcal{U}(3,n_2)} > 0$ (see Conjecture~\ref{conj:unique3}). However, general formulas for ~${\mathcal{U}(2,n_2)}$ and ~${\mathcal{U}(3,n_2)}$ remain to be derived (see Open~\ref{open:fillInMore}). Finally, using a brute force computation, we calculated the number of unique stable states for~${n_1=2}$ when $n_2\in[8]$ and ${n_1=3}$ when $n_2\in[7]$. The results are given in Table~\ref{tab:NumberOfUniqueStableStates}. 

\begin{table}
	\centering
	\begin{tabular}{c|rrrrrrrrrr}
		\diagbox{$n_1$}{$n_2$} & $1$ & $2$ & $3$ & $4$ & $5$ & $6$ & $7$ & $8$ & \ldots & $n_2$\\
		\hline
		$1$ & $1$ & $2$ & $6$ & $24$ & $120$ & $720$ & $5,040$ & $40,320$ & \ldots & $n_2!$\\
		$2$ &   & $12$ & $106$ & $1,108$ &  $12,826$ & $163,276$ & $2,274,592$ & $34,318,068$ & $\ldots$ & $>0$\\
		$3$ &   &   & $966$ & $1,484$ & $3,528$ & $8,176$ & $18,592$ & $???$& $\ldots$ & $???$\\
		$4$ &   &   &    & 0 &  0 & 0 & 0 & 0 & $\ldots$ & 0\\
		$\ge 5$ &   &   &    &    & 0 & 0 & 0 & 0 & $\ldots$ & 0\\
	\end{tabular}
 \vspace{.2 cm}
	\caption{
	    A partially filled chart giving $\mathcal{U}(n_1, n_2)$. for small $n_1$ and $n_2$. If Conjecture~\ref{conj:unique3} holds, then the question marks in row 3 must be strictly positive. 
	}
	\label{tab:NumberOfUniqueStableStates}
\end{table}

Despite calculating ${\mathcal{U}(2,n_2)}$ and~${\mathcal{U}(3,n_2)}$ for small values of $n_2$, the sequences~${\mathcal{U}(2,n_2)}$ and~${\mathcal{U}(3,n_2)}$ remain mysterious, and do not match any existing OEIS entries. 
\begin{open}
\label{open:fillInMore}
    What are the sequences~${\mathcal{U}(2,n_2)}$ and~${\mathcal{U}(3,n_2)}$? Can these be expressed in a closed form or at least asymptotically?
\end{open}

While we do not know much about the sequence ${\mathcal{U}(3,n_2)}$, we are fairly confident that it is always positive. In particular, we conjecture the following. 
\begin{conjecture}\label{conj:unique3}
    For any $n_2 \in \Z_{\ge 1}$, we have $\mathcal U(3,n_2) >0$. 
\end{conjecture}
One way to prove Conjecture~\ref{conj:unique3} is to construct an explicit point set with a unique stable state for each $n_2$. We have a candidate construction. Let $\sigma_1 = \identityperm$, and define $\sigma_2$ in the following way:
    \[\sigma_2(j) =  \begin{cases}3j-2 & \text{for $1 \le j \le n_2$}\\
    3(2n_2 + 1 - j) & \text{for $n_2 < j \le 2n_2$}\\
    3(j-2n_2)-1 & \text{for $2n_2 < j \le 3n_2$}
   \end{cases}.\] 
The construction above gives a unique stable state for $n_2 \le 7$, and we conjecture that this generalizes to all $n_2$ (which would prove Conjecture~\ref{conj:unique3}).

Note that Table~\ref{tab:NumberOfUniqueStableStates} gives an almost complete picture of unique stable states for two-dimensional grids.
While Lemma~\ref{lem:nc1} and~\ref{lem:nc2} as well as Proposition~\ref{pro:PropagatingNon-Uniqueness} hold in arbitrary dimension~$d$, it is a priori unclear how many unique stable states there are for grids of higher dimension~$d\geq3$ with~$n_i\leq3$ for all but one index~$i$.
Using the aforementioned computer searches, we found 1,823,944 point sets that admit a unique stable state for three-dimensional grids of size~$2\times2\times 2$.
However, the computational search quickly becomes extremely expensive.
It is therefore left as future work to expand this investigation to the higher-dimensional setting.

As a final remark in this section, note that in the proof of Proposition~\ref{pro:U(4,4)=0}, we found that each point set of~16 points allows for at least two stable placements.
However, this lower bound on the minimum number of stable placements is not necessarily sharp.
Using a computer search, we checked about 7\% of the point sets on the~${4\times 4}$ grid. 
The minimum number of stable states that one point set had, was seven. 
This minimum number came up several times, for different point sets.
However, as we did not pick the point sets uniformly random, we cannot extrapolate how many point sets have exactly seven stable states.
Also, there might be other point sets with even lower numbers of stable states that we have not found yet.
Hence, we ask the following, open question.

\begin{open}
\label{open:MinNumStableStates}
	For any specific grid with no point sets that have a unique stable state (e.g., the~${4 \times 4}$ grid), what is the minimum number of stable states a point set has? And is it possible to classify the point sets for which this minimum is achieved?
\end{open}

Having phrased this open question, we turn to its opposite: What are the point sets with the highest number of stable states?
We will investigate this question in the following section.


\section{Point Sets with Maximally Many Stable States}

Fix~$d$ and~$(n_1,\dots,n_d)$ with~${N = \prod n_i}$.
Recall that~$\identityperm$ denotes the identity permutation.
We will write~$P_{\identityperm}$ as shorthand for~$P_{(\identityperm, \dots, \identityperm)}$ (the point set associated with~$d$ copies of the same permutation).
We can think of a state of~$P_{\identityperm}$ as an assignment of a single number to each position on the grid by identifying each point~$(a,\ldots,a)$, $a\in[N]$ with the number~$a$. 
The state is stable if and only if these numbers are ordered in every Cartesian direction when traversing the grid.

Because of the reduction to numbers, $P_{\identityperm}$ is one-dimensional.
Hence, one might suspect that it would have a large number of stable states. 
In fact, we will show in Corollary~\ref{cor:identityMaximum} that when~${d=2}$, the number of stable states of~$P_{\identityperm}$ is maximal among all point sets. 
We conjecture that this is true for all~$d$ (see Conjecture~\ref{conj:identityworst}).

In order to show that the set of stable states of~$P_{\identityperm}$ is maximally large among all sets of stable states of other point sets~$P_{\bm\sigma}$ for any collection~$\bm\sigma$, consider the following algorithm.
It defines a map~${\varphi_{\bm\sigma}: \calS(P_{\bm\sigma}) \to \calS(P_{\identityperm})}$ that transforms every stable state of~$P_{\bm\sigma}$ into a stable state of~$P_{\identityperm}$.

\begin{algorithm}
	\caption{Stable State to Identity Stable State for $d=2$} 
	\label{alg:MaptoIdentity}
	\begin{algorithmic}[1]
		\Procedure{$\varphi_{\bm\sigma}: \calS(P_{\bm\sigma}) \to \calS(P_\I)$}{}
		\State We begin with a stable state $M_\pi(P_{\bm\sigma})$, and can assume that $\sigma_1 = \identityperm$. 
		\State Sort the values in each column by their first coordinate.
		\State Replace $\sigma_2$ with $\identityperm$. 
		\EndProcedure
	\end{algorithmic}
\end{algorithm}

\begin{example}
    To illustrate Algorithm~\ref{alg:MaptoIdentity}, consider, e.g., the~${3\times 3}$ stable state of~$P_{\bm\sigma}$ with~${\sigma_2=\langle 9,4,3,8,5,2,7,6,1\rangle}$ given in Figure~\ref{fig:NonEquivalenceXYBin}.
    Sorting the columns by the first coordinate brings the top element of the first and third column to the bottom while leaving the second column untouched. 
    Then, replacing the second coordinate by the first, i.e., reducing to one coordinate, yields a stable state of~$P_{\identityperm}$.
    See Figure~\ref{fig:IllustrationOfVarphi} for an illustration of this.
    Note how both columns and rows are ordered increasingly from bottom to top and left to right, respectively.
    
    \begin{figure}
        \centering
            \begin{tikzpicture}[scale=0.6]
                \draw[step=2.0] (0,0) grid (6,6);
    		\draw (1,1) node[scale=0.9] {(2,4)}; 
    		\draw (1,3) node[scale=0.9] {(4,8)}; 
    		\draw (1,5) node[scale=0.9] {(1,9)}; 
    		\draw (3,1) node[scale=0.9] {(3,3)}; 
    		\draw (3,3) node[scale=0.9] {(5,5)}; 
    		\draw (3,5) node[scale=0.9] {(7,7)}; 
    		\draw (5,1) node[scale=0.9] {(9,1)}; 
    		\draw (5,3) node[scale=0.9] {(6,2)}; 
    		\draw (5,5) node[scale=0.9] {(8,6)}; 
                \draw (8,3) node[scale=0.9] {$\stackrel{\substack{\text{Sort columns},\\\text{by the first}\\\text{coordinate}}}{\longrightarrow}$};
    	\end{tikzpicture}
            \begin{tikzpicture}[scale=0.6]
                \draw[step=2.0] (0,0) grid (6,6);
    		\draw (1,1) node[scale=0.9] {(1,9)}; 
    		\draw (1,3) node[scale=0.9] {(2,4)}; 
    		\draw (1,5) node[scale=0.9] {(4,8)}; 
    		\draw (3,1) node[scale=0.9] {(3,3)}; 
    		\draw (3,3) node[scale=0.9] {(5,5)}; 
    		\draw (3,5) node[scale=0.9] {(7,7)}; 
    		\draw (5,1) node[scale=0.9] {(6,2)}; 
    		\draw (5,3) node[scale=0.9] {(8,6)}; 
    		\draw (5,5) node[scale=0.9] {(9,1)}; 
                \draw (8.5,3) node {$\stackrel{\substack{\text{Replace second}\\\text{coordinate by first}\\\text{and reduce to it.}}}{\longrightarrow}$};
    	\end{tikzpicture}
            \begin{tikzpicture}[scale=0.6]
                \draw[step=2.0] (0,0) grid (6,6);
    		\draw (1,1) node {$1$};
    		\draw (1,3) node {$2$};
    		\draw (1,5) node {$4$};
    		\draw (3,1) node {$3$};
    		\draw (3,3) node {$5$};
    		\draw (3,5) node {$7$};
    		\draw (5,1) node {$6$};
    		\draw (5,3) node {$8$};
    		\draw (5,5) node {$9$};
            \end{tikzpicture}
        \caption{Illustrating the procedure given by Algorithm~\ref{alg:MaptoIdentity} when applied to a stable state from Figure~\ref{fig:NonEquivalenceXYBin}.}
        \label{fig:IllustrationOfVarphi}
    \end{figure}
\end{example}

The function~$\varphi_{\bm\sigma}$ from Algorithm~\ref{alg:MaptoIdentity} provides a way to relate the stable states of any point set~$P_{\bm\sigma}$ to those of~$P_{\identityperm}$.
The following theorem establishes that~$\varphi_{\bm\sigma}$ is injective and well-defined, i.e., that it indeed maps to~$\calS(P_{\identityperm})$.

\begin{theorem}
\label{thm:identitymost}
    The function~$\varphi_{\bm\sigma}$ defined in Algorithm~\ref{alg:MaptoIdentity} is injective and always maps to a stable state in~$\calS(P_{\identityperm})$.
\end{theorem}

\begin{proof}
    We first show that~$\varphi_{\bm\sigma}$ is injective. 
    Notice that~$\varphi_{\bm\sigma}$ preserves the first coordinate of the points in each column. 
    This means that for placement maps~$\pi$ and~$\pi'$, if~$\varphi_{\bm\sigma}(M_{\pi}(P_{\bm\sigma})) = \varphi_{\bm\sigma}(M_{\pi'}(P_{\bm\sigma}))$, then the states $M_{\pi}(P_{\bm\sigma})$ and $M_{\pi'}(P_{\bm\sigma})$ must have the same points in each column. 
    However, there is at most one way to arrange the points within each column to give a stable state. 
    Thus, we must have $\pi = \pi'$.  

    What is left to show is that~$\varphi_{\bm\sigma}$ maps stable states to stable states. 
    Let~${M_{\pi'}(P_{\identityperm}) = \varphi_{\bm\sigma}(M_\pi(P_{\bm\sigma}))}$. 
    By construction, the columns of~$M_{\pi'}(P_{\identityperm})$ must be sorted, so we just have to check the rows. 
    Recall that we can express~$M_{\pi'}(P_{\identityperm})$ by writing a single number in each box of our grid. 
    Assume for the sake of contradiction that there are two columns of~$M_{\pi'}(P_{\identityperm})$ with their respective entries as 
	\begin{align*}
		M_{\pi'}(P_{\identityperm})=
		\begin{array}{|c|c|c|c|c|}
			\hline
			\ldots & p_{n_2} & \ldots & q_{n_2} & \ldots\\
			\hline
			\vdots & \vdots & \iddots & \vdots & \vdots \\
			\hline
			\ldots & p_i & \ldots & q_i & \ldots\\
			\hline
			\vdots & \vdots & \iddots & \vdots & \vdots \\
			\hline
			\ldots & p_1 & \ldots & q_1 & \ldots\\
			\hline
		\end{array}
	\end{align*}
    such that there is some~${i\in[{n_2}]}$ with~${p_i>q_i}$. 
    Since the columns are sorted increasingly, we have~${p_i<p_{i+1}<\ldots<p_{n_2}}$. 
    Because~$\varphi_{\bm\sigma}$ did not change the position of the first coordinate between columns, for each~$p_j$, there is some~$q_k$ such that~$p_j$ and~$q_k$ were placed in the same row in~${M_\pi(P_{\bm\sigma})}$. 
    Since~${M_\pi}(P_{\bm\sigma})$ is stable, it follows that~${p_j<q_k}$. 
    Thus, there are~${{n_2}-i+1}$ values~$q_k$ larger than~$p_i$. 
    However, if~${q_i < p_i}$, then all larger values have to be above~$q_i$ in the column as the column is sorted increasingly. 
    There are only~${{n_2}-i}$ places in the column left to store the~${{n_2}-i+1}$ larger values, which gives the desired contradiction.
\end{proof}

By this theorem, we have an injective map from the set of stable states of any point set~$P_{\bm\sigma}$ to the set of stable states of~$P_{\identityperm}$.
As these sets are finite, the following corollary follows immediately.

\begin{corollary}
	\label{cor:identityMaximum}
	For $d=2$, the size of the set~$\calS(P_{\bm\sigma})$ is maximized when~${{\bm\sigma} = (\identityperm,\ldots,\identityperm)}$. 
\end{corollary}

There are several natural ways to extend the definition of~$\varphi_{\bm\sigma}$ to all~${d \in \mathbb N}$, but we were not able to find a map satisfying both of the properties of Theorem~\ref{thm:identitymost}. 
Still, we take it as a basis to formulate the following conjecture.

\begin{conjecture}\label{conj:identityworst}
Corollary~\ref{cor:identityMaximum} holds for all $d \in \mathbb N$. 
\end{conjecture}

For the case of~${d=2}$, stable states for the identity permutation have been extensively studied in the the context of \emph{standard Young tableaux}, which we will briefly discuss below.

\begin{definition}
    Let~$\lambda_1, \lambda_2, \dots, \lambda_k$ be a non-increasing collection of positive integers. 
    The \emph{Young diagram} $\lambda = (\lambda_1, \lambda_2, \dots, \lambda_k)$ is the diagram formed by $k$ rows of boxes where the $i$-{th} row contains $\lambda_i$ boxes. 
\end{definition}

To match our grid conventions, will write these diagrams using \emph{French notation}, where the rows are left-justified and are written upward from widest to narrowest (see Figure~\ref{fig:YoungTableau}). 
We will be interested in the Young diagram formed by~$n_1$ rows each with~$n_2$ boxes each (i.e., the Young diagram~$\lambda=(n_2^{n_1})$, where this notation is short for ${\lambda=(\underbrace{n_2, \ldots, n_2}_{n_1\text{ times}})}$).

\begin{definition}
Given a Young diagram with~$N$ boxes, a \emph{standard Young tableau (on the alphabet $\{1,\ldots,N\}$)} is an assignment of the integers~$1$ to~$N$ to the boxes of the Young diagram such that the entries increase along rows and up columns.
\end{definition}

\begin{example}
    In Figure~\ref{fig:YoungTableau}, we consider the Young diagram for~${N=4}$ and the partition~${\lambda=(3,1)}$ with two different associated Young tableaux. The tableau on the left is standard while the tableau on the right is not. A straightforward combinatorial calculation shows that there are 24 fillings of this Young diagram, 3 of which are standard. 
\end{example}

\begin{figure}
	\centering
	\begin{tikzpicture}[scale=0.6]
	\draw (0,0) grid (3,1);
	\draw (0,1) grid (1,2);
	\draw (0.5,0.5) node {1};
	\draw (1.5,0.5) node {2};
	\draw (2.5,0.5) node {4};
	\draw (0.5,1.5) node {3};
	\end{tikzpicture}
	\hspace{2cm}
	\begin{tikzpicture}[scale=0.6]
	\draw (0,0) grid (3,1);
	\draw (0,1) grid (1,2);
	\draw (0.5,0.5) node {4};
	\draw (1.5,0.5) node {1};
	\draw (2.5,0.5) node {3};
	\draw (0.5,1.5) node {2};
	\end{tikzpicture}
	\caption{Two Young tableaux corresponding to the partition~${\lambda=(3,1)}$ of~${N=4}$. The left tableau is standard, while the right tableau is not. }
	\label{fig:YoungTableau}
\end{figure}

By definition, an assignment of the numbers~$1$ through~$N$ to the Young diagram~$\lambda=(n_2^{n_1})$ produces a standard Young tableau if and only if this assignment corresponds to a stable state~$M_\pi(P_{\bm\sigma})$ for~${\bm{\sigma} = (\identityperm,\identityperm)}$.
The number of standard Young tableaux for a given partition of~$N$ is given by the \emph{hook length formula}, which was first proven by Frame, Robinson, and Thrall~\cite{frame1954hook}. 
See~\cite[Section 10]{sagan2001symmetric} for more information and history about this formula. 
The following result is immediate from the hook length formula. 

\begin{proposition}\label{prop:youngnum}
    The number of standard Young tableaux associated with $\lambda = (n_2^{n_1})$ is precisely
        \begin{align*}
            f^{n_1\times n_2}:=\frac{(n_1 n_2)!}{\prod_{i=1}^{n_1}\prod_{j=1}^{n_2}(n_1 + n_2 - i - j + 1)}.
	\end{align*}
\end{proposition}

\begin{corollary}\label{cor:YoungNumIdentityMax}
When $d=2$, the value of $|\calS(P_{\identityperm})|$ is given by the formula from Proposition~\ref{prop:youngnum}. This is also an upper bound for $|\calS(P_{\bm \sigma})|$ among all pairs of permutations $\bm \sigma$. 
\end{corollary}
\begin{proof}
    The first claim is immediate from Proposition~\ref{prop:youngnum}. The second claim is immediate from Corollary~\ref{cor:identityMaximum}. 
\end{proof}

Recall that with Corollary~\ref{cor:U(n_1,...,n_d)=0}, we have shown that the assumption made by Malheiros and Walter~\cite{malheiros2015simple} to prove the equivalent of Corollary~\ref{cor:timeOptimality} was not justified.
We gave a complete proof of Theorem~\ref{the:lowerBound} above and thereby proved Corollary~\ref{cor:timeOptimality}.
However, we obtain an alternative proof of Corollary~\ref{cor:timeOptimalityD=2} by using Corollary~\ref{cor:YoungNumIdentityMax}. 
\begin{proof}
Recall that the proof of Theorem~\ref{the:lowerBound} built a decision-tree and bounded the minimal depth that this tree needs to be traversed to by any comparison-based algorithm~$\mathcal{A}$ in order to construct a placement~$\pi$ such that~$M_\pi(P_{\bm\sigma})$ is stable.
From Corollary~\ref{cor:YoungNumIdentityMax}, we know that any such decision-tree has to be traversed at least until the subtree has at most~$|\mathcal{S}(P_{\identityperm})|$ leaves. 
That is, assuming~${n_1=n_2=n}$ for some~${n\in\mathbb{N}}$, the traversal has to take at least this many steps:
\begin{align*}
	  &\log(n^2!) - \log(f^{n\times n})\\
	= &\log\left(\frac{n^2!}{f^{n\times n}}\right)\\
	= & \log\left(\prod_{i=1}^n\prod_{j=1}^n (2n-i-j+1)\right)\\
	= & \sum_{i=1}^n\sum_{j=1}^n \log(2n-i-j+1)\\
	= & \log(2n-1)+\log(2n-2)+\ldots+\log(n)\\
	  & + \log(2n-2) + \log(2n-3) + \ldots + \log(n-1)\\
	  & \vdots\\
	  & + \log(n) + \log(n-1) + \ldots + \log(1)\\
	= & n\log(n) + \sum_{i=1}^{n-1}\left(i\log(i)+i\log(2n-i)\right)\\
	\geq & n\log(n) + \frac{1}{2}\sum_{i=n/2}^{n-1} \frac{n}{2}\log(n)\\
	= & \Omega(n^2\log(n)).
\end{align*}
Thereby, we have obtained the statement as made by~\cite{malheiros2015simple} in two different ways, for the two-dimensional case.
\end{proof}

In the setting of~${d=2}$, the map~$\varphi_{\bm\sigma}$, as constructed by Algorithm~\ref{alg:MaptoIdentity}, gives rise to an algorithm for the enumeration of all stable states of a given point set. 
This enumeration algorithm is significantly faster than the naive brute force approach. 
Consider the procedure given in Algorithm~\ref{alg:EnumerateStableStates}.

\begin{algorithm}
	\caption{Enumeration of Stable States} 
	\label{alg:EnumerateStableStates}
	\begin{algorithmic}[1]
		\Procedure{Enumerate Stable States}{$P_{\bm\sigma}$}
		\State Enumerate all stable states~${M_\pi(\mathcal{I})}$ of the identity~$\mathcal{I}$
		\For{$M_\pi(\mathcal{I})$ stable}
			\State Compute~$\varphi_{\bm\sigma}^{-1}(M_\pi(\mathcal{I}))$
			\State Store it, if stable
		\EndFor
		\EndProcedure
	\end{algorithmic}
\end{algorithm}

Given this algorithm, enumeration of all stable states~$M_\pi(\identityperm)$ takes time~$\mathcal{O}(f^{n_1\times n_2})$.
Execution of Algorithm~\ref{alg:MaptoIdentity} then takes time~$\mathcal{O}(n_1 n_2\log(n_2))$.
Furthermore, if $n_2>n_1$, then we can reverse the role of $n_1$ and $n_2$ without changing the result. Thus, the total time to execute Algorithm~\ref{alg:EnumerateStableStates} amounts to
\begin{align}
    f^{n_1\times n_2} n_1 n_2\log(\min(n_1,n_2)).
\end{align}

For a concrete example of the benefits of this approach, consider the case where ${n_1=n_2=4}$. Given a point set~$P_{\bm \sigma}$, a naive method for enumerating stable states is to check all $16! \approx 2.09\cdot10^{13}$ possible placement maps. Using Algorithm~\ref{alg:EnumerateStableStates}, it suffices to compute just~$870,912,000$, a reduction by five orders of magnitude. One question that we expect this reduction to simplify is Open~\ref{open:MinNumStableStates}.

We will now turn to the final contribution of this section.
In Corollary~\ref{cor:identityMaximum}, we proved that~${\bm{\sigma}=(\identityperm,\identityperm)}$ maximizes the number of stable states for~${d=2}$.
However, this gives a mere example for a collection of permutations~$\sigma$ that realizes this maximum.
Below, we will characterize a family of point sets~$P_{\bm\sigma}$ that maximize~$|\mathcal{S}(P_{\bm\sigma})|$. We conjecture that this is a complete characterization apart from the trivial case where $n_1 = 1$ or $n_2 = 1$, and one exceptional case where $n_1 = n_2 = 2$ (see Conjecture~\ref{conj:fibon}). 

First, we note the following relationship about the cardinality of two related sets of stable states, which holds in arbitrary dimension~$d$.

\begin{lemma}\label{lem:revperm}
    Let~$\bm{\sigma}=(\sigma_1,\dots,\sigma_d)$ be a collection of permutations defining the point set~$P_{\bm\sigma}$. 
    For any~$k$, if we replace~$\sigma_k$ with the reverse permutation~$\overline{\sigma}_k$, we do not change the number of stable states.
\end{lemma}
\begin{proof}
    We obtain a bijection between~$\mathcal{S}(P_{(\sigma_1,\ldots,\sigma_k,\ldots,\sigma_d)})$ and~$\mathcal{S}(P_{(\sigma_1,\ldots,\overline{\sigma}_k,\ldots,\sigma_d)})$ as follows.
    Consider a stable state~$M_\pi(P_{\bm\sigma})$.
    By reverting the order of all points along all~$k$-strips, we obtain a stable state~$M_{\pi'}(P_{(\sigma_1,\ldots,\overline{\sigma}_k,\ldots,\sigma_d)})$.
    Note that this operation is bijective.
\end{proof}

Now, we focus on the two-dimensional case.
Here, we can prove the following relationship between the number of stable states of two related point sets.

\begin{lemma}\label{lem:swapadj}
    Let $\bm \sigma = (\identityperm,\sigma_2)$. 
    Suppose for some~${i\in [N]}$ that~${|\sigma_2(i) -  \sigma_2(i+1)| = 1}$. 
    Define~$\bm \sigma'$ to be equivalent to $\bm \sigma$, except with~${\sigma'_2(i) = \sigma_2(i+1)}$ and~${\sigma'_2(i+1) = \sigma_2(i)}$. 
    Then, the point set~$P_{\bm \sigma'}$ has the same number of stable states as~$P_{\bm \sigma}$. 
\end{lemma}
\begin{proof}
    We will define a function~$\psi$ which gives a permutation on the placement maps. 
    Given a placement map~$\pi$, let~${\bm a = (a_1,a_2)}$ be the position given by~$\pi^{-1}(i)$ and let~${\bm b = (b_1,b_2)}$ be the position given by~$\pi^{-1}(i+1)$. 
    For any position~${\bm c \in [n_1] \times [n_2]}$, let
    \begin{align*}
        \psi\circ \pi (\bm c) = 
        \begin{cases} 
            \bm b &\text{if $\bm c = \bm a$ and $a_1 = b_1$,}\\
            \bm a &\text{if $\bm c = \bm b$ and $a_1 = b_1$,}\\
            \bm c &\text{otherwise.}
        \end{cases}
    \end{align*}
    In other words,~$\psi$ switches the position of~${(i,\sigma_2(i))}$ and~${(i+1,\sigma_2(i+1))}$ when these points are in the same column, and leaves everything the same otherwise. 

    We claim that~$M_\pi(P_{\bm\sigma})$ is stable if and only if~$M_{\psi\circ\pi}(P_{\bm\sigma'})$ is stable. The result then follows from the fact that $\psi$ is invertible. 
    
    To prove the claim, there are two possibilities to consider. First, if $a_1 = b_1$, then the only difference between $M_\pi(P_{\bm\sigma})$ and $M_{\psi\circ\pi}(P_{\bm\sigma'})$ is that the first entry of $(i,\sigma(i))$ and $(i+1,\sigma(i+1))$ are swapped. By the assumption that $a_1 = b_1$, these points must be in the same column and thus must be in different rows. Combining this with the fact that there are no numbers between $i$ and $i+1$ implies that this change does not impact stability. 
    
    The other case to consider is where $a_1 \not= b_1$. This time, the only difference between $M_\pi(P_{\bm\sigma})$ and $M_{\psi\circ\pi}(P_{\bm\sigma'})$ is that the second entry of $(i,\sigma(i))$ and $(i+1,\sigma(i+1))$ are swapped. By the condition that~${|\sigma_2(i) -  \sigma_2(i+1)| = 1}$, there can be no values between $\sigma_2(i)$ and  $\sigma_2(i+1)$. The claim follows from combining this fact with the fact that the swapped entries are not in the same column. 
\end{proof}

\begin{lemma}\label{lem:whichfuncs} A point set $\bm \sigma = (\identityperm,\sigma_2)$ can be obtained from $(\identityperm,\identityperm)$ by repeated applications of Lemma~\ref{lem:swapadj} if and only if $|\sigma_2(i) - i|\le 1$ for all $i\in[N]$. 
\end{lemma}
\begin{proof} The forward direction is straightforward to prove by induction. If $\sigma_2(1) = 1$, then we can ignore the first entry. Alternatively, if $\sigma_2(1) = 2$, then we must have $\sigma_2(2) = 1$ and we can ignore the first two entries (after a single swap for $i=1$). 

For the reverse direction, it suffices to show that if $|\sigma_2(j) - j|\le 1$ for all $j\in [N]$, then $|\sigma'_2(j) - j|\le 1$ for all $j\in [N]$ (using the notation from Lemma~\ref{lem:swapadj}). Let $i$ be the index such that $\sigma'_2(i) = \sigma_2(i+1)$ and $\sigma'_2(i+1) = \sigma_2(i)$. For $k \not\in \{i,i+1\}$, it is immediate that $|\sigma'_2(k) - k| = |\sigma_2(k) - k|\le 1$. Thus, we must prove that $|\sigma'_2(i) - i| \le 1$ and $|\sigma'_2(i+1) - (i+1)| \le 1$. 

By the conditions of Lemma~\ref{lem:swapadj}, the only way to arrive at a contradiction would be if $\sigma_2(i) = i+1$ and $\sigma_2(i+1) = i+2$ or if $\sigma_2(i) = i-1$ and $\sigma_2(i+1) = i$. However, if $\sigma_2(i) = i+1$ and $\sigma_2(i+1) = i+2$, then for all $k > i+1$, we must have $\sigma_2(k) \in [i+3,N]$. This is impossible by the pigeonhole principle, and the other case is analogous. 
\end{proof}

Lemmas~\ref{lem:revperm} and~\ref{lem:whichfuncs} can also be used to produce enumerative results. 

\begin{definition}
    Let $(n_1,\ldots,n_d)\in\mathbb{N}^d$, then we write $\mathcal{M}(n_1,\dots,n_d)$ for the number of point sets with the maximal number of stable states on the~${(n_1 \times \dots \times n_d)}$ grid. 
\end{definition}

Let $\fibon_k$ denote the $k^{th}$ Fibonacci number (with $\fibon_1 =\fibon_2 =1$, and $\fibon_k =\fibon_{k-1} + \fibon_{k-2}$). It is known that for any $N$, there are exactly $\fibon_{N+1}$ functions satisfying Lemma~\ref{lem:whichfuncs}. This follows from a simple recursion after splitting into the case where $\sigma(1) = 1$ and the case where $\sigma(1) = 2$. By Lemma~\ref{lem:swapadj}, this gives $\fibon_{N+1}$ point sets with the same number of stable states as the identity. Furthermore, as long as $N >2$, we obtain $\fibon_{N+1}$ more point sets with this property by applying Lemma~\ref{lem:revperm}. Thus, we have the following corollary. 

\begin{corollary}\label{cor:fibonMin} 
    For $(n_1,n_2)\in \mathbb N^2$ with $n_1n_2>2$, we have $\mathcal M(n_1,n_2) \ge 2\fibon_{n_1n_2+1}$, where $\fibon_k$ denotes the $k^{th}$ Fibonacci number. 
\end{corollary}

Based on computer data and heuristic arguments, we suspect that this bound is almost always tight. In particular, we conjecture the following. 

\begin{conjecture}\label{conj:fibon} 
    For $(n_1,n_2)\in \mathbb N^2$ with $n_1>1$, $n_2>1$, and $n_1n_2>4$, we have $\mathcal M(n_1,n_2) = 2\fibon_{n_1n_2+1}$,
\end{conjecture}

Note that if $n_1 =1$ or $n_2 = 1$, then all $N!$ point sets have a unique stable state. When $n_1 = n_2 = 2$, Corollary~\ref{cor:fibonMin} gives $10$ point sets with the maximum of $2$ stable states. However, there are also $2$ more point sets with this property corresponding to the permutations $\sigma_2 = (2,4,1,3)$ and $\sigma_2 = (3,1,2,4)$. 

With the help of a computer, we verified that $\mathcal M(2,3) = 26 = 2\fibon_{7}$, and we suspect that since Conjecture~\ref{conj:fibon} holds for this grid, it should hold for larger grids as well.

\section{Conclusion}

In this article, we have investigated the neighborhood grid data structure~\cite{joselli2009neighborhood,joselli2015neighborhood} from a combinatorial viewpoint.
This perspective allowed us to show that an ordered arrangement of points in the grid, a stable state, can be constructed in polynomial time for any grid shape in any dimension, see Theorem~\ref{the:StableStateExistence}).
In terms of run time of the construction of stable states, we proved a lower bound, see Theorem~\ref{the:lowerBound}, and proved time optimality of the construction algorithm for the case where all side lengths of the grid are equal, see Corollary~\ref{cor:timeOptimality}.
Previous work did provide corresponding arguments, which, however, rested on a faulty uniqueness assumption.
We disprove this assumption by showing that all point sets placed in grids of size~$4\times 4$ or larger admit to at least two stable states, see Propositions~\ref{pro:U(4,4)=0} and~\ref{pro:PropagatingNon-Uniqueness}.
However, we can also bound the number of stable states that a point set admits to from above, see Proposition~\ref{prop:youngnum}.

In the two-dimensional case, for one grid side being of length~${n_1=2}$, we were able to prove that there always exist certain point sets that have a unique stable state.
However, it is an open question to derive a closed-form expression of the number of these stable states, see Open~\ref{open:fillInMore}.
Similarly, for~${n_1=3}$, an exact formula is unknown and it remains an open question whether there are unique stable states for any~$n_2$.
We do have a candidate point set that only allows for one unique placement in all instances we checked, but we were not able to prove this generally.
This candidate motivates our Conjecture~\ref{conj:unique3} that for~$n_3$ there are also always unique stable states for any~$n_2$.
On the other hand, given a larger grid, with side lengths~${n_1,n_2\geq4}$, it remains unclear what the minimum number of stable states is among all respective point sets, see Open~\ref{open:MinNumStableStates}. 
Furthermore, while several of our results hold for arbitrary dimension~$d$, the investigation of unique stable states in higher-dimensional configurations is left as future work.

We also proved that the maximum number of stable states is achieved by the identity permutation in the two-dimensional case, and suspect that this is true for arbitrary dimension, see Conjecture~\ref{conj:identityworst}. After this, we gave a partial classification of point sets which also achieve the maximum set by the identity permutation, which we suspect is a full classification unless $n_1 = 1$, $n_2 = 1$, or $n_1 = n_2 = 2$, see Conjecture~\ref{conj:fibon}. This conjecture suggests an elegant formula for the number of such point sets related to the Fibonacci sequence. 

Other questions rise in the investigation of the neighborhood grid data structure. 
These include, for instance, the neighborhood quality achieved by the grid.
As neighborhood queries are only answered approximately, it remains unclear how close the reported neighbors are.
Furthermore, in an application context, points representing, e.g., moving particles, will not move for large distances in simulations.
Thus, it is probably advised to use an adaptive strategy based on the grid in the previous iteration and not to re-build the grid in every iteration of the simulation.
However, these questions are out scope for the combinatorial consideration presented here and will have to be tackled in a separate work.

\section*{Acknowledgments}

This material is based upon work supported by the National Science Foundation under Grant No.\@ DMS-1439786 and the Alfred P. Sloan Foundation award G-2019-11406 while the author was in residence at the Institute for Computational and Experimental Research in Mathematics in Providence, RI, during the Illustrating Mathematics program.
This research was partially funded by the Deutsche Forschungsgemeinschaft (DFG, German Research Foundation) -- 455095046.
We acknowledge the support by the DFG SFB/TR 109 `Discretization in Geometry and Dynamics', ECMath, BMS, and the German National Academic Foundation. 
The authors would like to thank Shagnik Das for helpful discussions regarding the proof of Theorem~\ref{the:lowerBound} for the case of~${n_i=n}$ and~${d=2}$, i.e., the statement that is now Corollary~\ref{cor:timeOptimality}.
Furthermore, the authors would like to thank Mathijs Molenaar and Mark van de Ruit for creating a C++ implementation to experimentally investigate stable states. Finally, we would like to thank Per Alexandersson and Sam Hopkins for answering some questions we had about higher-dimensional generalizations of Young tableaux.

\bibliographystyle{plain}
\bibliography{sources}

\appendix

\section*{Appendix}

In the following appendix, we collect results and experiments regarding the neighborhood grid that are not necessarily of combinatorial nature.

\section{Iterative Parallelized Procedure}
\label{sec:IterativeParallelizedProcedure}

A benefit of the neighborhood grid data structure not discussed so far is the straight forward parallelization of an algorithm creating a stable state. In this section, we discuss the parallelization of a building algorithm for stable states and investigate its runtime.

\subsection{Parallelization of Stable State Creation}
\label{subsec:ParallelizationOfStableStateCreation}

The idea of iterative and parallel creation of stable states is discussed in~\cite{malheiros2015simple} at great detail and we are only going to state the basic ideas and results here. However, some questions are not investigated, which we are going to tackle in this article. The general scheme for a two-dimensional grid is quite simple and given in Algorithm~\ref{alg:IterativeParallelizedSorting}.

\begin{algorithm}
	\caption{Iterative Parallelized Sorting of a Neighborhood Grid} \label{alg:IterativeParallelizedSorting}
	\begin{algorithmic}[1]
		\Procedure{Iterative Parallelized Sorting}{Grid $G$}
		\While{$\neg$sorted($G$)}
		\State sort all rows of $G$ in parallel
		\State sort all columns of $G$ in parallel
		\EndWhile
		\EndProcedure
	\end{algorithmic}
\end{algorithm}

An immediate question following this algorithm concerns its convergence. Does it converge and if so, is the final state always stable? Malheiros and Walter did not answer this question in~\cite{malheiros2015simple}, nor did Joselli et al.\@ in the original publications~\cite{joselli2009neighborhood,joselli2015neighborhood}, where the idea of Algorithm~\ref{alg:IterativeParallelizedSorting} is also used. We will prove the convergence here.

\begin{theorem}[U. Reitebuch, M. Skrodzki]
	\label{the:ParallelConvergence}
	Given any point set~${P=\{p_1,\ldots,p_N\mid p_i\in\mathbb{R}^2\}}$, placing the points in a matrix~$M$ and running Algorithm~\ref{alg:IterativeParallelizedSorting} on the matrix, the algorithm converges and yields a stable state of~$M$.
\end{theorem}
\begin{proof}
	Given a matrix~$M$, consider the following expression:
	\begin{align}
	\label{equ:Energy}
	E(M)=\sum_{i,j=1}^n i\cdot a_{ij}+j\cdot b_{ij},
	\end{align}
	for each sorting step of Algorithm~\ref{alg:IterativeParallelizedSorting}, this expression grows, but it can at most attain~$N!$ many different values. Hence, the algorithm converges.
	
	If the matrix is not in a stable state, i.e.\@ there is a row or column violating the stable state, sorting this row or column lets expression~(\ref{equ:Energy}) grow and resolves the conflict in the given row or column, possibly creating a new conflict in another row or column. Therefore, a local maximum of this expression is equivalent to a stable state in the matrix.
\end{proof}

Note that the energy functional~(\ref{equ:Energy}) can easily be extended to higher-dimensional settings. Given a corresponding version of Algorithm~\ref{alg:IterativeParallelizedSorting} for the higher-dimensional case, Theorem~\ref{the:ParallelConvergence} holds true for arbitrary dimensional point sets.

\begin{figure}
	\includegraphics[width=1.0\textwidth]{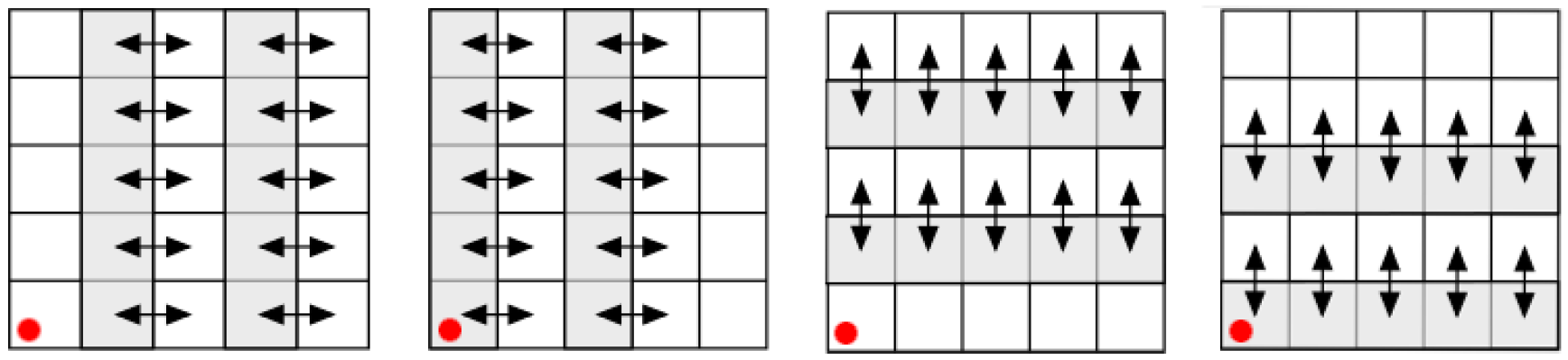}
	\caption{Step-based odd-even sort on a two-dimensional neighborhood grid, figure from \cite{malheiros2015simple}.}
	\label{fig:OddEvenSort}
\end{figure}

In their work~\cite{malheiros2015simple}, Malheiros and Walter investigate a slight variation of Algorithm~\ref{alg:IterativeParallelizedSorting}. Namely, they do not perform a full sorting of a row or column, but rather consider one step of the odd-even sort algorithm of Habermann, see~\cite{habermann1972parallel}. Such step performs an exchange between all those cells in odd columns and their respective right neighboring cells, if this pair of cells violates the stable state conditions of Definition~\ref{def:StableState}. Then, all cells in even columns and their respective right neighboring cells are compared and exchanged if necessary. The same is performed on all odd rows and even rows, yielding a four step mechanism, see Figure~\ref{fig:OddEvenSort}.

Note that the argument in the proof of Theorem~\ref{the:ParallelConvergence} holds true also for this algorithm. Therefore, it also converges to a stable state. Question remains how fast the algorithm works asymptotically. Note that, defying intuition, elements can cycle using this procedure. An example is given in Figure~\ref{fig:IllustrateCycling}, where the algorithm depicted in Figure~\ref{fig:OddEvenSort} is used.

\begin{figure}
	\begin{tikzpicture}
		\draw (0,0.5) rectangle (2,2.5);
		\draw (0,1.5) -- (2,1.5);
		\draw (1,0.5) -- (1,2.5);
		\node at (0.5,1) {$(0,3)$};
		\node at (0.5,2) {$(3,1)$};
		\node at (1.5,1) {$(1,0)$};
		\node[red] at (1.5,2) {$(2,2)$};
	\end{tikzpicture}
	$\stackrel{x}{\rightarrow}$
	\begin{tikzpicture}
		\draw (0,0.5) rectangle (2,2.5);
		\draw (0,1.5) -- (2,1.5);
		\draw (1,0.5) -- (1,2.5);
		\node at (0.5,1) {$(0,3)$};
		\node[red] at (0.5,2) {$(2,2)$};
		\node at (1.5,1) {$(1,0)$};
		\node at (1.5,2) {$(3,1)$};
	\end{tikzpicture}
	$\stackrel{y}{\rightarrow}$
	\begin{tikzpicture}
		\draw (0,0.5) rectangle (2,2.5);
		\draw (0,1.5) -- (2,1.5);
		\draw (1,0.5) -- (1,2.5);
		\node[red] at (0.5,1) {$(2,2)$};
		\node at (0.5,2) {$(0,3)$};
		\node at (1.5,1) {$(1,0)$};
		\node at (1.5,2) {$(3,1)$};
	\end{tikzpicture}
	$\stackrel{x}{\rightarrow}$
	\begin{tikzpicture}
		\draw (0,0.5) rectangle (2,2.5);
		\draw (0,1.5) -- (2,1.5);
		\draw (1,0.5) -- (1,2.5);
		\node at (0.5,1) {$(1,0)$};
		\node at (0.5,2) {$(0,3)$};
		\node[red] at (1.5,1) {$(2,2)$};
		\node at (1.5,2) {$(3,1)$};
	\end{tikzpicture}
	$\stackrel{y}{\rightarrow}$
	\begin{tikzpicture}
		\draw (0,0.5) rectangle (2,2.5);
		\draw (0,1.5) -- (2,1.5);
		\draw (1,0.5) -- (1,2.5);
		\node at (0.5,1) {$(1,0)$};
		\node at (0.5,2) {$(0,3)$};
		\node at (1.5,1) {$(3,1)$};
		\node[red] at (1.5,2) {$(2,2)$};
	\end{tikzpicture}
	\caption{Performing an odd-even sort alternating on all the rows and columns causes the red element~${(2,2)}$ to cycle through the matrix.}
	\label{fig:IllustrateCycling}
\end{figure}

Given the possibility of cycling elements, the theoretical asymptotic bounds of the parallel iterative algorithm remain unclear. Only experimental results are available, as presented in~\cite{malheiros2015simple}. 

\begin{open}
\label{open:WorstCaseParallelRuntime}
	Given the step-wise odd-even algorithm depicted in Figure~\ref{fig:OddEvenSort}, what is its parallelized worst-case runtime aside from the upper bound of~${\mathcal{O}(N!)}$ as established in the proof of Theorem~\ref{the:ParallelConvergence}?
\end{open}
First experiments and intuition lead us to state the following conjecture:
\begin{conjecture}
	The step-wise odd-even algorithm depicted in Figure~\ref{fig:OddEvenSort} runs like Bubble-Sort in worst case time of~${\mathcal{O}(N^2)}$, after parallelization in~${\mathcal{O}(N)}$.
\end{conjecture}

\subsection{Comparison to k-d Trees and Adaptive Sorting Algorithms for Fast Modifications}
\label{subsec:AdaptiveSortingAlgorithmsForFastModifications}

The most wide-spread data structure for neighborhood computation, k-d~trees, popular for its expected nearest neighbor lookup time of~${\mathcal{O}(\log(n))}$ (see~\cite{friedman1977algorithm}), suffers from a severe problem. Namely, if the underlying point set is slightly altered, the k-d~tree might become unbalanced. Although there are some heuristics how to modify k-d~trees when adding or deleting points, at some stage the k-d~tree has to be rebuild, which is costly. The authors of~\cite{gaede1998multidimensional} conclude to this end:
\begin{quote}
	\textit{The adaptive k-d~tree is a rather static structure; it is obviously difficult to keep the tree balanced in the presence of frequent insertions and deletions.}
\end{quote}
In contrast, the neighborhood grid is built only with sorting algorithms. For these, adaptive algorithms are available that benefit from a sorted set into which a small number of records is to be inserted, cf.~\cite{petersson1995framework}. Thus, altering the point set underlying the neighborhood grid can be performed faster than rebuilding a k-d~tree. An exact investigation of this relation is left for future research, see open question~\ref{open:AdaptiveSortingBenefit}.
\begin{open}
	\label{open:AdaptiveSortingBenefit}
	Is utilizing an adaptive sorting algorithm on a modified stable state faster than rebuilding a whole new stable state?
\end{open}

Note that the algorithm presented in Theorem~\ref{the:StableStateExistence} needs to sort the given points. When utilizing~${N/2}$ processors, sorting can be performed in~$\log(N)$ time, see~\cite{ajtai1983sorting}. Therefore, the presented algorithm can be parallelized to run in~${\mathcal{O}(\log(N))}$. This particular approach is of rather theoretical relevance, as the constants in~\cite{ajtai1983sorting} are comparably large. However, other work is devoted to finding more practical parallelizations, see~\cite{amato1996comparison}. Also, it makes for a significant speed-up compared to the algorithm depicted in Figure~\ref{fig:OddEvenSort}.

Compare this to building a k-d~tree in parallel. A straight-forward parallelization would be as follows: In each step~$i$, we have to sort~$i$ sets of~$n^2/2^i$ points in the dimension with largest spread, which takes~$\log(n^2)-i$ time for each of the~$\log(n^2)$ levels of the tree, resulting in an upper bound for the total building time of~$\Omega(\log^2(n))$. However, this only holds for a straight forward parallelization. As a k-d~tree can be used for sorting by placing all points along one dimension, by~\cite{leighton1984tight}, it has a minimum build time of~$\mathcal{O}(\log(n))$. To the best of our knowledge, it is unclear whether the gap between this lower bound and the upper bound induced by the straight forward parallelization can be closed.

Therefore, the neighborhood grid can be built slightly faster when compared with the straight forward parallelized k-d~tree, but only gives estimated answers, while the k-d~tree provides exact neighbor relations.

\section{Quality of Neighborhood Approximation}
\label{sec:QualityOfNeighborhoodApproximation}

As stated above, the neighborhood estimates given by the neighborhood grid data structure are not necessarily precise. In this section, we will investigate the quality of the neighborhood approximation.

\subsection{Single Point Neighbor}
\label{subsec:SinglePointNeighbor}

A first question to answer in this subsection concerns the distance of two geometrical neighbors from~$P$ in the stable state of~$M$. In the following, we will create a point set~$P$ with points~${p,q\in P}$ that are respective nearest neighbors to each other, but that lie on the exact opposite sites of~$M$.

Consider the following points
\begin{align*}
	&{p=(0,0)}, {q=(1,0)},\\
	&{p_i=(0,2+i/n)}, && {i=1,\ldots,n-1},\\
	&{q_i=(1,-2-i/n)}, && {i=1,\ldots,n-1},\\ 
	&{r_{i,j}=(1-1/i,2+j/n)}, && {i=2,\ldots,n-1}, {j=1,\ldots,n}.
\end{align*} 
This yields a point set~$P$ with~$n^2$ points for~${n\geq2}$, see Figure~\ref{fig:SingPointNeighborPointSet}. Given these points, the following matrix is in a stable state:
\begin{align*}
M(P)=
\begin{array}{|c|c|c|c|c|}
\hline
p_{n-1} & r_{2,n} & \ldots & r_{n-1,n} & q\\
\hline
\vdots & \vdots & \iddots & \vdots & q_1\\
\hline
p_1 & r_{2,2} & \ldots & r_{n-1,2} & \vdots\\
\hline
p & r_{2,1} & \ldots & r_{n-1,n} & q_{n-1}\\
\hline
\end{array}.
\end{align*}

\begin{figure}
	\centering
	\begin{tikzpicture}
		\draw[gray!50] (-2,-3) grid (2,4);
		\draw (0,0) circle (1.0);
		\draw (1,0) circle (1.0);
		\filldraw[black] (0,0) circle (2pt)  node[above left] {$p$}; 
		\filldraw[black] (1,0) circle (2pt) node[above right] {$q$}; 
		\filldraw[black] (0,2.25) circle (2pt)  node[left] {$p_1$}; 
		\filldraw[black] (0,2.5) circle (2pt)  node[left] {$p_2$}; 
		\filldraw[black] (0,2.75) circle (2pt)  node[left] {$p_3$}; 
		\filldraw[black] (1,-2.25) circle (2pt)  node[right] {$q_1$}; 
		\filldraw[black] (1,-2.5) circle (2pt)  node[right] {$q_2$}; 
		\filldraw[black] (1,-2.75) circle (2pt)  node[right] {$q_3$}; 
		\filldraw[black] (0.5,2.25) circle (2pt); 
		\filldraw[black] (0.5,2.5) circle (2pt); 
		\filldraw[black] (0.5,2.75) circle (2pt); 
		\filldraw[black] (0.5,3) circle (2pt); 
		\filldraw[black] (0.66,2.25) circle (2pt); 
		\filldraw[black] (0.66,2.5) circle (2pt); 
		\filldraw[black] (0.66,2.75) circle (2pt); 
		\filldraw[black] (0.66,3) circle (2pt) node[above] {$r_{i,j}$}; 
	\end{tikzpicture}
	\caption{Point set~$P$ as given in Section~\ref{subsec:SinglePointNeighbor}.}
	\label{fig:SingPointNeighborPointSet}
\end{figure}

Note that the nearest neighbor to~$p$ and~$q$ in~$P$ is~$q$ and~$p$ respectively. However, in~$M(P)$, these points lie in the opposing corners of the matrix. That is, in order to find the geometrically closest neighbor to~$p$ in~$M(P)$, the~$n$-ring around~$p$ has to be checked. In other words, all points have to be checked, which takes~${\mathcal{O}(n^2)}$ instead of~${\theta(\log(n))}$ as in k-d~trees.

\subsection{All Point Nearest Neighbors}
\label{subsec:AllPointNearestNeighbors}

In the previous example we saw that for a single point, its unique nearest neighbor can be arbitrarily far away in the neighborhood grid. However, when considering all points, how is the overall estimate? Here, we provide an example, where no point has its corresponding neighbor within its one-ring in the neighborhood grid.

For~$n\in\mathbb{N}$,~$2\mid n$, consider the following points
\begin{align*}
	&p_{i,j} = (i,j), && q_{i,j}=(i+n,j+0.5),\\
	&r_{i,j}=(i+0.5,j+n), && 	s_{i,j}=(i+n+0.5,j+n+0.5),\\ &i,j\in\{0,\ldots,\frac{n}{2}-1\}.
\end{align*}	
This yields a point set~$P$ with~$n^2$ points, see Figure~\ref{fig:AllPointNeighborPointSet}. Given these points, the following matrix is in a stable state:
\begin{align*}
M(P)=
\begin{array}{|c|c|c|c|c|c|c|}
\hline
\textcolor{red}{q_{0,\frac{n}{2}-1}} & \textcolor{ForestGreen}{s_{0,\frac{n}{2}-1}} & \textcolor{red}{q_{1,\frac{n}{2}-1}} & \textcolor{ForestGreen}{s_{1,\frac{n}{2}-1}} & \ldots & \textcolor{red}{q_{\frac{n}{2}-1,\frac{n}{2}-1}} & \textcolor{ForestGreen}{s_{\frac{n}{2}-1,\frac{n}{2}-1}}\\
\hline
p_{0,\frac{n}{2}-1} & \textcolor{blue}{r_{0,\frac{n}{2}-1}} & p_{1,\frac{n}{2}-1} & \textcolor{blue}{r_{1,\frac{n}{2}-1}} & \ldots & p_{\frac{n}{2}-1,\frac{n}{2}-1} & \textcolor{blue}{r_{\frac{n}{2}-1,\frac{n}{2}-1}}\\
\hline
\vdots & \vdots & \vdots & \vdots & \iddots &\vdots&\vdots\\
\hline
\textcolor{red}{q_{0,1}} & \textcolor{ForestGreen}{s_{0,1}} & \textcolor{red}{q_{1,1}} & \textcolor{ForestGreen}{s_{1,1}} & \ldots & \textcolor{red}{q_{\frac{n}{2}-1,1}} & \textcolor{ForestGreen}{s_{\frac{n}{2}-1,1}}\\
\hline
p_{0,1} & \textcolor{blue}{r_{0,1}} & p_{1,1} & \textcolor{blue}{r_{1,1}} & \ldots & p_{\frac{n}{2}-1,1} & \textcolor{blue}{r_{\frac{n}{2}-1,1}}\\
\hline
\textcolor{red}{q_{0,0}} & \textcolor{ForestGreen}{s_{0,0}} & \textcolor{red}{q_{1,0}} & \textcolor{ForestGreen}{s_{1,0}} & \ldots & \textcolor{red}{q_{\frac{n}{2}-1,0}} & \textcolor{ForestGreen}{s_{\frac{n}{2}-1,0}}\\
\hline
p_{0,0} & \textcolor{blue}{r_{0,0}} & p_{1,0} & \textcolor{blue}{r_{1,0}} & \ldots & p_{\frac{n}{2}-1,0} & \textcolor{blue}{r_{\frac{n}{2}-1,0}}\\
\hline
\end{array}.
\end{align*}

\begin{figure}
	\centering
	\begin{tikzpicture}
	\draw[gray!50] (-1,-1) grid (6,6);
	\filldraw[black] (0,0) circle (2pt)  node[below left] {$p_{0,0}$}; 
	\filldraw[black] (1,0) circle (2pt)  node[below left] {$p_{1,0}$}; 
	\filldraw[black] (0,1) circle (2pt)  node[below left] {$p_{0,1}$}; 
	\filldraw[black] (1,1) circle (2pt)  node[below left] {$p_{1,1}$}; 

	\filldraw[red] (4,0.5) circle (2pt)  node[below left] {$q_{0,0}$}; 
	\filldraw[red] (5,0.5) circle (2pt)  node[below left] {$q_{1,0}$}; 
	\filldraw[red] (4,1.5) circle (2pt)  node[below left] {$q_{0,1}$}; 
	\filldraw[red] (5,1.5) circle (2pt)  node[below left] {$q_{1,1}$}; 
	
	\filldraw[blue] (0.5,4) circle (2pt)  node[below left] {$r_{0,0}$}; 
	\filldraw[blue] (1.5,4) circle (2pt)  node[below left] {$r_{1,0}$}; 
	\filldraw[blue] (0.5,5) circle (2pt)  node[below left] {$r_{0,1}$}; 
	\filldraw[blue] (1.5,5) circle (2pt)  node[below left] {$r_{1,1}$}; 
	
	\filldraw[ForestGreen] (4.5,4.5) circle (2pt)  node[below left] {$s_{0,0}$}; 
	\filldraw[ForestGreen] (5.5,4.5) circle (2pt)  node[below left] {$s_{1,0}$}; 
	\filldraw[ForestGreen] (4.5,5.5) circle (2pt)  node[below left] {$s_{0,1}$}; 
	\filldraw[ForestGreen] (5.5,5.5) circle (2pt)  node[below left] {$s_{1,1}$}; 
	\end{tikzpicture}
	\caption{Point set~$P$ as given in Section~\ref{subsec:AllPointNearestNeighbors} for~$n=4$.}
	\label{fig:AllPointNeighborPointSet}
\end{figure}

Note that the nearest neighbor to~${p_{i,j}}$ is some~${p_{k,\ell}}$, the nearest neighbor to~${q_{i,j}}$ is some~${q_{k,\ell}}$, etc. However, none of the points has its corresponding neighbor in its one-ring. Expanding this scheme by adding more four-point subsets, it is easily achievable to create point sets with immediate stable states where no respective nearest neighbor is in the~${\frac{n}{4}}$ ring of all points.

\subsection{Average Nearest Neighbor Quality}

We have given two constructions to create worst case behavior in the neighborhood estimates. However, these constructions do not account for the average estimate to be expected for a random point set. Therefore, we pose the following questions:

\begin{open}
	\label{open:QualityOfNeighborhoodEstimate}
	How good is the average neighborhood estimate of the neighborhood grid?
	How does the building procedure affect the neighborhood quality?
\end{open}

\begin{figure}
	\centering
	\includegraphics[width=1.0\textwidth]{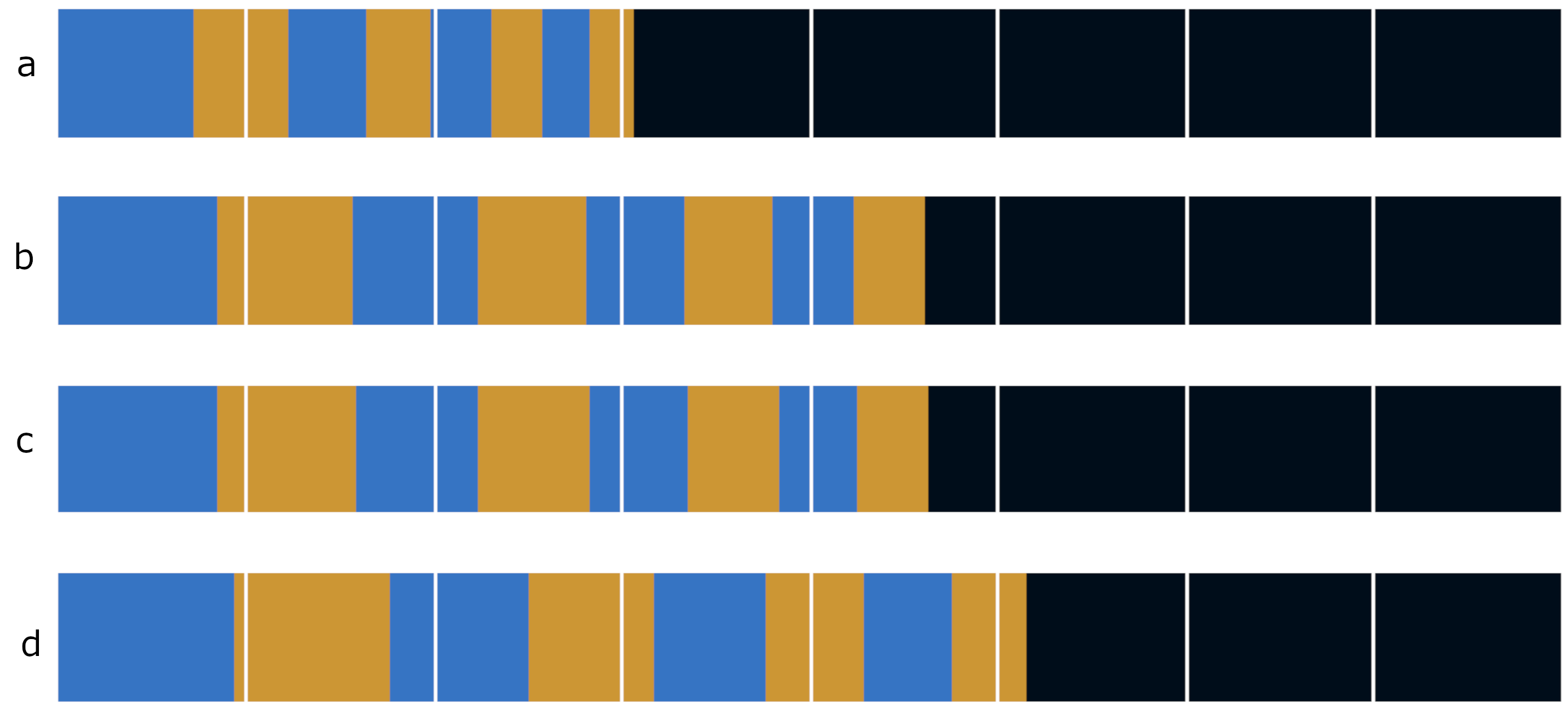}
	\caption{From top to bottom: Percentage of neighbors correctly reporting their 1st, 2nd,~$\ldots$, 8th nearest neighbors when using (a) direct sorting following Theorem~\ref{the:StableStateExistence}, (b) performing complete odd-even-sort of~\cite{habermann1972parallel} on all rows and columns alternating, (c) performing one step odd-even-sort of~\cite{malheiros2015simple} on all rows and columns alternating, (d) iteratively exchanging the pair of points which provide largest grow in Energy~(\ref{equ:Energy}. From left to right, the boxes indicate whether every point knows its nearest, its two nearest,~$\ldots$, its eight nearest neighbors among its one-ring of eight points. The colored bars indicate to what percentage the respective boxes are filled. Note that the direct method of Theorem~\ref{the:StableStateExistence} runs fastest among these methods but provides the worst results. Exchanging two points for maximum grow of
	Energy~(\ref{equ:Energy}) is the most time-consuming procedure but provides the best results.}
	\label{fig:ExperimentalResults}
\end{figure}

Concerning both questions, we present experimental results here. We generated~$70$ point sets with~$1,764$ points each randomly in~${[0,1]^2\subset \mathbb{R}^2}$. They are randomly placed into the grid. Then, each of the following four methods is applied in order to build a stable state on the grid:
\begin{enumerate}
	\item[(a)] The grid is sorted directly following the algorithm from Theorem~\ref{the:StableStateExistence}.
	\item[(b)] All rows are sorted, then all columns are sorted, the procedure is iterated until convergence.
	\item[(c)] The step-based odd-even sort of~\cite{habermann1972parallel} is performed.
	\item[(d)] We iteratively exchange those points~$p_i,p_j$ in the grid such that the grow in Energy~(\ref{equ:Energy}) is maximal.
\end{enumerate}
As the output of procedures (b)-(d) heavily depends on the initial placement of the points in the grid, we apply the procedure on~$30$ different random placements of each point set. Finally, we count how many points in the grid have their nearest neighbor, their two nearest neighbors, $\ldots$, and their eight nearest neighbors within their one-ring. The results are shown in Figure~\ref{fig:ExperimentalResults}. Note that the most time-consuming methods---exchanging by maximum growth of Energy~(\ref{equ:Energy})---gives the best results while the direct sorting following Theorem~\ref{the:StableStateExistence} performs worst despite being the fastest sorting method.

\end{document}